\documentclass[amsfonts, 12pt]{article}
\usepackage{amssymb}

\newtheorem{theorem}{Theorem}[section]

\newtheorem{example}[theorem]{Example}
\newtheorem{proposition}[theorem]{Proposition}
\newtheorem{lemma}[theorem]{Lemma}

\newtheorem{remark}[theorem]{Remark}

\def\cB{\mathcal{B}}
\def\cC{\mathcal{C}}

\def\cG{\mathcal{G}}

\def\cF{\mathcal{F}}
\def\cL{\mathcal{L}}

\def\cN{\mathcal{N}}
\def\cO{\mathcal{O}}
\def\cP{\mathcal{P}}

\def\bE{\mathbb{E}}
\def\bR{\mathbb{R}}

\def\ve{\varepsilon}

\begin{document}

\title{SPDEs with $\alpha$-stable L\'evy noise: \\ a random field approach}

\author{Raluca M. Balan\footnote{Department of Mathematics and Statistics, University of Ottawa,
585 King Edward Avenue, Ottawa, ON, K1N 6N5, Canada. E-mail
address: rbalan@uottawa.ca} \footnote{Research supported by a
grant from the Natural Sciences and Engineering Research Council
of Canada.}}

\date{November 25, 2013}
\maketitle

\begin{abstract}
\noindent This article is dedicated to the study of an SPDE of the form $$Lu(t,x)=\sigma(u(t,x))\dot{Z}(t,x) \quad t>0, x \in \cO$$ with zero initial conditions and Dirichlet boundary conditions, where $\sigma$ is a Lipschitz function, $L$ is a second-order pseudo-differential operator, $\cO$ is a bounded domain in $\bR^d$, and $\dot{Z}$ is an $\alpha$-stable L\'evy noise with $\alpha \in (0,2)$, $\alpha\not=1$ and possibly non-symmetric tails. To give a meaning to the concept of solution, we develop a theory of stochastic integration with respect to $Z$, by generalizing the method of \cite{GM83} to higher dimensions and non-symmetric tails. The idea is to first solve the equation with ``truncated'' noise $\dot{Z}_{K}$ (obtained by removing from $Z$ the jumps which exceed a fixed value $K$), yielding a solution $u_{K}$, and then show that the solutions $u_L,L>K$ coincide
on the event $t \leq \tau_{K}$, for some stopping times $\tau_K \uparrow \infty$ a.s. A similar idea was used in \cite{PZ07} in the setting of Hilbert-space valued processes. A major step is to show that the stochastic integral with respect to $Z_{K}$
satisfies a $p$-th moment inequality, for $p \in (\alpha,1)$ if $\alpha<1$, and $p \in (\alpha,2)$ if $\alpha>1$. This inequality plays the same role as the Burkholder-Davis-Gundy inequality in the theory of integration with respect to continuous martingales.
\end{abstract}



\section{Introduction}

Modeling phenomena which evolve in time or space-time and are subject to random perturbations is a fundamental problem in stochastic analysis. When these perturbations are known to exhibit an extreme behavior, as seen frequently in finance or environmental studies, a model relying on the Gaussian distribution is not appropriate. A suitable alternative could be a model based on a heavy-tailed distribution, like the stable distribution. In such a model, these perturbations are allowed to have extreme values with a probability which is significantly higher than in a Gaussian-based model.

In the present article, we introduce precisely such a model, given rigorously by a stochastic partial differential equation (SPDE) driven by a noise term which has a stable distribution over any space-time region, and has independent values over disjoint space-time regions (i.e. it is a L\'evy noise). More precisely, we consider the SPDE:
\begin{equation}
\label{eq-noiseZ}
Lu(t,x)=\sigma(u(t,x))\dot{Z}(t,x), \quad t>0, x \in \cO
\end{equation}
with zero initial conditions and Dirichlet boundary conditions, where $\sigma$ is a Lipschitz function, $L$ is a second-order pseudo-differential operator on a bounded domain $\cO \subset \bR^d$, and $\dot{Z}(t,x)=\frac{\partial^{d+1}Z}{\partial t \partial x_1 \ldots \partial x_d}$ is the formal derivative of an $\alpha$-stable L\'evy noise with $\alpha \in (0,2)$, $\alpha \not=1$. The goal is to find sufficient conditions on the fundamental solution $G(t,x,y)$ of the equation $Lu=0$ on $\bR_{+} \times \cO$, which will ensure the existence of a mild solution of equation (\ref{eq-noiseZ}). We say that a predictable process $u=\{u(t,x);t \geq 0,x \in \cO\}$ is a {\bf mild solution} of (\ref{eq-noiseZ}) if for any $t>0,x \in \cO$,
\begin{equation}
\label{def-sol}
u(t,x)=\int_0^t \int_{\cO}G(t-s,x,y)\sigma(u(s,y))Z(ds,dy) \quad \mbox{a.s.}
\end{equation}
We assume that $G(t,x,y)$ is a function in $t$, which excludes from our analysis the case of the wave equation with $d \geq 3$.

To explain the connections with other works, we describe briefly the construction of the noise (the details are given in Section \ref{def-noise-section} below). This construction is similar to that of a classical $\alpha$-stable L\'evy process, and is based on a Poisson random measure (PRM) $N$ on $\bR_{+} \times \bR^d \times (\bR \verb2\2 \{0\})$ of intensity $dt dx \nu_{\alpha}(dz)$, where
\begin{equation}
\label{def-nu-alpha}\nu_{\alpha}(dz)=[p \alpha z^{-\alpha-1}1_{(0,\infty)}(z)+q \alpha (-z)^{-\alpha-1}1_{(-\infty,0)}(z)]dz
\end{equation}
for some $p,q \geq 0$ with $p+q=1$.
More precisely, for any set $B \in \cB_{b}(\bR_{+} \times \bR^d)$,
\begin{equation}
\label{represent-Z}
Z(B)=\int_{B \times \{|z| \leq 1\}}z \widehat{N}(ds,dx,dz)+\int_{B \times \{|z| > 1\}}zN(ds,dx,dz)-\mu|B|,
\end{equation}
where $\widehat{N}(B \times \cdot)=N(B \times \cdot)-|B| \nu_{\alpha}(\cdot)$ is the compensated process and $\mu$ is a constant (specified by Lemma \ref{Y(B)param} below). Here, $\cB_{b}(\bR_{+} \times \bR^d)$ is the class of bounded Borel sets in $\bR_{+} \times \bR^d$ and $|B|$ is the Lebesque measure of $B$.

As the term on the right-hand side of (\ref{def-sol}) is a stochastic integral with respect to $Z$, such an integral should be constructed first.
Our construction of the integral is an extension to random fields of the construction provided by Gin\'e and Marcus in \cite{GM83} in the case of an $\alpha$-stable L\'evy process $\{Z(t)\}_{t \in [0,1]}$. Unlike these authors, we do not assume that the measure $\nu_{\alpha}$ is symmetric. 

Since any L\'evy noise is related to a PRM, in a broad sense, one could say that this problem originates in It\^o's papers \cite{ito44} and \cite{ito46} regarding the stochastic integral with respect to a Poisson noise.
SPDEs driven by a compensated PRM were considered for the first time in \cite{kallianpur-et-al94}, using the approach based on Hilbert-space-valued solutions. This study was motivated by an application to neurophysiology leading to the cable equation.
In the case of the heat equation, a similar problem was considered in \cite{AWZ98}, \cite{bie98} and \cite{applebaum-wu00} using the approach based on random-field solutions.
One of the results of \cite{bie98} shows that the heat equation:
$$\frac{\partial u}{\partial t}(t,x)=\frac{1}{2}\Delta u(t,x)+\int_{U}f(t,x,u(t,x);z)\widehat{N}(t,x,dz)+g(t,x,u(t,x))$$
has a unique solution in the space of predictable processes $u$ satisfying \linebreak
$\sup_{(t,x) \in [0,T] \times \bR^d}E|u(t,x)|^p<\infty$, for any $p \in (1+2/d,2]$.
In this equation, $\widehat{N}$ is the compensated process corresponding to
a PRM $N$ on $\bR_{+} \times \bR^d \times U$ of intensity $dtdx \nu(dz)$, for an arbitrary $\sigma$-finite measure space $(U,\cB(U), \nu)$ with measure $\nu$ satisfying $\int_{U}|z|^p \nu(dz)<\infty$. Because of this later condition, this result cannot be used in our case with $U=\bR \verb2\2 \{0\}$ and $\nu=\nu_{\alpha}$. For similar reasons, the results of \cite{applebaum-wu00} also do not cover the case of an $\alpha$-stable  noise.
However, in the case $\alpha>1$, we will be able to exploit successfully some ideas of \cite{bie98} for treating the equation with ``truncated'' noise $Z_K$, obtained by removing from $Z$ the jumps exceeding a value $K$ (see Section \ref{tr-noise-alpha-greater-section} below).

The heat equation with the same type of noise as in the present article was examined in \cite{mueller98} and \cite{mytnik02} in the cases $\alpha<1$, respectively $\alpha>1$, 
assuming that the noise has only positive jumps (i.e. $q=0$). The methods used by these authors are different from those presented here, since they investigate the more difficult case of a non-Lipschitz function $\sigma(u)=u^{\delta}$ with $\delta>0$.
In \cite{mueller98}, Mueller removes the atoms of $Z$ of mass smaller than $2^{-n}$ and solves the equation driven by the noise obtained in this way; here we remove the atoms of $Z$ of mass larger than $K$ and solve the resulting equation.
In \cite{mytnik02}, Mytnik uses a martingale problem approach and gives the existence of a pair $(u,Z)$ which satisfies the equation (the so-called ``weak solution''), whereas in the present article we obtain the existence of a solution $u$ for a {\em given} noise $Z$ (the so-called ``strong solution''). In particular, when $\alpha>1$ and $\delta=1/\alpha$, the existence of a ``weak solution'' of the heat equation with $\alpha$-stable L\'evy noise is obtained in \cite{mytnik02} under the condition
\begin{equation}
\label{heat-cond}
\alpha<1+\frac{2}{d}
\end{equation}
which we encounter here as well. It is interesting to note that (\ref{heat-cond}) is the necessary and sufficient condition for the existence of the density of the super-Brownian motion with ``$\alpha-1$''-stable branching (see \cite{dawson92}). Reference \cite{mueller-mytnik-stan06} examines the heat equation with multiplicative noise (i.e. $\sigma(u)=u$), driven by an $\alpha$-stable L\'evy noise $Z$ which does not depend on time.

To conclude the literature review, we should point out that there are many references related to stochastic differential equations with $\alpha$-stable L\'evy noise, using the approach based on Hilbert-space valued solutions. We refer the reader to Section 12.5 of the monograph \cite{PZ07}, and to \cite{PZ06}, \cite{AMR09}, \cite{marinelli-rockner10},  \cite{priola-zabczyk11} for a sample of relevant references. See also the survey article \cite{oksendal08} for an approach based on the white noise theory for L\'evy processes.

This article is organized as follows.
\begin{itemize}
\item In Section \ref{def-noise-section}, we review the construction of the $\alpha$-stable L\'evy noise $Z$, and we show that this can be viewed as an independently scattered random measure with jointly $\alpha$-stable distributions.
\item In Section \ref{linear-eq-section}, we consider the linear equation (\ref{eq-noiseZ}) (with $\sigma(u)=1$) and we identify the necessary and sufficient condition for the existence of the solution. This condition is verified in the case of some examples.
\item Section \ref{stoch-int-section} contains the construction of the stochastic integral with respect to the $\alpha$-stable noise $Z$, for $\alpha \in (0,2)$. The main effort is dedicated to proving a maximal inequality for the tail of the integral process, when the integrand is a simple process. This extends the construction of \cite{GM83} to the case random fields and non-symmetric measure $\nu_{\alpha}$.
\item In Section \ref{truncated-noise-section}, we introduce the process $Z_K$ obtained by removing from $Z$ the jumps exceeding a fixed value $K$, and we develop a theory of integration with respect to this process. For this, we need to treat separately the cases $\alpha<1$ and $\alpha>1$. In both cases, we obtain a $p$-th moment inequality for the integral process for $p \in (\alpha,1)$ if $\alpha<1$, and $p \in (\alpha,2)$ if $\alpha>1$. This inequality plays the same role as the Burkholder-Davis-Gundy inequality in the theory of integration with respect to continuous martingales.
\item In Section \ref{nonlinear-eq-section} we prove the main result about the existence of the mild solution of equation (\ref{eq-noiseZ}). For this, we first solve the equation with ``truncated'' noise $Z_K$ using a Picard iteration scheme, yielding a solution $u_K$. We then introduce a sequence $(\tau_K)_{K \geq 1}$ of stopping times with $\tau_K \uparrow \infty$ a.s. and we show that the solutions $u_L,L>K$ coincide on the event $t \leq \tau_K$. For the definition of the stopping times $\tau_K$, we need again to consider separately the cases $\alpha<1$ and $\alpha>1$.
\item Appendix \ref{appA-section} contains some results about the tail of a non-symmetric stable random variable, and the tail of an infinite sum of random variables. 
    Appendix \ref{appB-section} gives an estimate for the Green function associated to the fractional power of the Laplacian. Appendix \ref{appC-section} gives a local property of the stochastic integral with respect to $Z$ (or $Z_K$).
\end{itemize}

\section{Definition of the noise}
\label{def-noise-section}

In this section we review the construction of the $\alpha$-stable L\'evy noise on $\bR_{+} \times \bR^d$ and investigate some of its properties.

Let $N=\sum_{i \geq 1}\delta_{(T_i,X_i,Z_i)}$ be a Poisson random measure on $\bR_{+} \times \bR^d \times (\bR \verb2\2 \{0\})$, defined on a probability space $(\Omega,\cF,P)$, with intensity measure $dtdx\nu_{\alpha}(dz)$, where $\nu_{\alpha}$ is given by (\ref{def-nu-alpha}). Let $(\ve_j)_{j \geq 0}$ be a sequence of positive real numbers such that $\ve_j \to 0$ as $j \to \infty$ and $1=\ve_0>\ve_1>\ve_2>\ldots$. Let
$$\Gamma_j=\{z \in \bR; \ve_{j}<|z| \leq \ve_{j-1}\}, \ j\geq 1 \quad \mbox{and} \quad \Gamma_0=\{z \in \bR; |z|>1\}.$$

For any set $B \in \cB_b(\bR_{+} \times \bR^d)$, we define
$$L_j(B)=\int_{B \times \Gamma_j}z N(dt,dx,dz)=\sum_{(T_i,X_i) \in B}Z_i 1_{\{Z_i \in \Gamma_j\}}, \quad j \geq 0.$$

\begin{remark}
{\rm The variable $L_0(B)$ is finite since the sum above contains finitely many terms.  To see this, we note that $E[N(B \times \Gamma_0)]=|B|\nu_{\alpha}(\Gamma_0)<\infty$, and hence $N(B \times \Gamma_0)={\rm card}\{i \geq 1; (T_i,X_i,Z_i) \in B \times \Gamma_0\}<\infty$.}
\end{remark}

For any $j \geq 0$, the variable $L_j(B)$ has a compound Poisson distribution with jump intensity measure $|B|\cdot \nu_{\alpha}|_{\Gamma_j}$, i.e.
\begin{equation}
\label{ch-funct-Lj(B)}
E[e^{iu L_j(B)}]=\exp\left\{|B|\int_{\Gamma_j}(e^{iuz}-1)\nu_{\alpha}(dz) \right\}, \quad u \in \bR.
\end{equation}

\noindent It follows that
$E(L_j(B))= |B| \int_{\Gamma_j} z \nu_{\alpha}(dz)$ and ${\rm Var}(L_j(B))= |B| \int_{\Gamma_j} z^2 \nu_{\alpha}(dz)$  for any $j \geq 0$.
Hence ${\rm Var}(L_j(B))<\infty$ for any $j \geq 1$ and ${\rm Var}(L_0(B))=\infty$. If $\alpha>1$, then $E(L_0(B))$ is finite.
Define
\begin{equation}
\label{def-Y(B)}
Y(B)=\sum_{j \geq 1}[L_j(B)-E(L_j(B))]+L_0(B).
\end{equation}
This sum converges a.s. by Kolmogorov's criterion since $\{L_j(B)-E(L_j(B))\}_{j \geq 1}$ are independent zero-mean random variables with
$\sum_{j \geq 1}{\rm Var}(L_j(B))
<\infty$.

From (\ref{ch-funct-Lj(B)}) and (\ref{def-Y(B)}), it follows that $Y(B)$ is an infinitely divisible random variable with characteristic function:
\begin{equation}
\label{chf-Y-B}
E(e^{iuY(B)})=\exp \left\{|B|\int_{\bR}(e^{iuz}-1-iu z 1_{\{|z| \leq 1\}})\nu_{\alpha}(dz) \right\}, \quad \ u \in \bR.
\end{equation}

\noindent Hence $E(Y(B))= |B|\int_{\bR}z 1_{\{|z|>1\}}\nu_{\alpha}(dz)$ and
${\rm Var}(Y(B))= |B|\int_{\bR}z^2 \nu_{\alpha}(dz)$.

\begin{lemma}
\label{Y-indep-scattered}
The family $\{Y(B);B \in \cB_{b}(\bR_{+} \times \bR^d)\}$ defined by (\ref{def-Y(B)}) is an independently scattered random measure, i.e.\\
(a) for any disjoint sets $B_1, \ldots,B_n$ in $\cB_{b}(\bR_{+} \times \bR^d)$, $Y(B_1), \ldots,Y(B_n)$ are independent ;\\
(b) for any sequence $(B_n)_{n \geq 1}$ of disjoint sets in $\cB_{b}(\bR_{+} \times \bR^d)$ such that $\bigcup_{n \geq 1}B_n$ is bounded,
$Y(\bigcup_{n \geq 1} B_n)=\sum_{n \geq 1}Y(B_n)$ a.s.
\end{lemma}

\noindent {\bf Proof:} (a) Note that for any function $\varphi \in L^2(\bR_{+} \times \bR^d)$ with compact support $K$, we can define the random variable
$Y(\varphi)=\sum_{j \geq 1}[L_j(\varphi)-E(L_j(\varphi))]+L_0(\varphi)$
where
$L_j(\varphi)=\int_{K \times \Gamma_j}\varphi(t,x)zN(dt,dx,dz)$.
For any $u \in \bR$, we have:
\begin{equation}
\label{chf-Y-varphi}
E(e^{iu Y(\varphi)})=\exp \left\{\int_{\bR_{+} \times \bR^d \times \bR}(e^{iu z \varphi(t,x)}-1-iuz\varphi(t,x)1_{\{|z| \leq 1\}}) dtdx \nu_{\alpha}(dz)\right\}.
\end{equation}

For any disjoint sets $B_1, \ldots, B_n$ and for any $u_1,\ldots,u_n \in \bR$, we have:
\begin{eqnarray}
\nonumber
\lefteqn{E[\exp({i\sum_{k=1}^n u_k Y(B_k)})]=E[\exp(i Y(\sum_{k=1}^{n} u_k 1_{B_k}))]}\\
\nonumber
& =& \exp \left\{\int_{\bR_{+} \times \bR^d \times \bR} (e^{iz \sum_{k=1}^{n} u_k 1_{B_k}(t,x)}-1-iz1_{\{|z| \leq 1\}} \sum_{k=1}^{n} u_k 1_{B_k}(t,x) )dt dx \nu_{\alpha}(dz)\right\}\\
\label{chf-calcul}
&=& \exp \left\{ \sum_{k=1}^{n}|B_k| \int_{\bR}(e^{iu_kz}-1-i u_k z 1_{\{|z| \leq 1\}})\nu_{\alpha}(dz)\right\}\\
\nonumber
&=& \prod_{k=1}^{n}E[\exp(iu_k Y(B_k))],
\end{eqnarray}
using (\ref{chf-Y-varphi}) with $\varphi=\sum_{k=1}^{n}u_k 1_{B_k}$ for the second equality, and (\ref{ch-funct-Lj(B)}) for the last equality. This proves that $Y(B_1), \ldots,Y(B_n)$ are independent.

(b) Let $S_n=\sum_{k=1}^{n}Y(B_k)$ and $S=Y(B)$, where $B=\bigcup_{n \geq 1}B_n$. By L\'evy's equivalence theorem, $(S_n)_{n \geq 1}$ converges a.s. if and only if it converges in distribution. By (\ref{chf-calcul}), with $u_i=u$ for all $i=1, \ldots,k$, we have:
$$E(e^{iuS_n})=\exp\left\{ |\bigcup_{k=1}^{n}B_k| \int_{\bR}(e^{iuz}-1-i u z 1_{\{|z| \leq 1\}})\nu_{\alpha}(dz)\right\}.$$
This clearly converges to $E(e^{iuS})=\exp\left\{|B| \int_{\bR}(e^{iuz}-1-i u z 1_{\{|z| \leq 1\}})\nu_{\alpha}(dz)\right\}$, and
hence $(S_n)_{n \geq 1}$ converges in distribution to $S$. $\Box$

\vspace{3mm}

Recall that a random variable $X$ has an {\em $\alpha$-stable distribution} with parameters $\alpha \in (0,2), \sigma \in [0,\infty),  \beta \in [-1,1],\mu \in \bR$ if for any $u \in \bR$,
\begin{eqnarray*}
E(e^{iuX})&=&\exp\left\{-|u|^{\alpha} \sigma^\alpha\left(1-i{\rm sgn}(u) \beta \tan \frac{\pi \alpha}{2}\right) +iu\mu\right\} \quad \mbox{if} \quad \alpha \not=1, \ \mbox{or}\\
E(e^{iuX})&=&\exp\left\{-|u| \sigma \left(1+i{\rm sgn}(u) \beta \frac{2}{\pi} \ln|u|\right) +iu\mu\right\} \quad \mbox{if} \quad \alpha =1
\end{eqnarray*}
(see Definition 1.1.6 of \cite{ST94}). We denote this distribution by $S_{\alpha}(\sigma,\beta,\mu)$.

\begin{lemma}
\label{Y(B)param}
$Y(B)$ has a $S_{\alpha}(\sigma|B|^{1/\alpha},\beta,\mu|B|)$ distribution with $\beta=p-q$,
$$\sigma^{\alpha}=\int_0^{\infty}\frac{\sin x}{x^{\alpha}}dx=\left\{
\begin{array}{ll}
\frac{\Gamma(2-\alpha)}{1-\alpha}\cos\frac{\pi \alpha}{2} & \mbox{if $\alpha \not=1$} \\
\frac{\pi}{2} & \mbox{if $\alpha=1$}
\end{array} \right., \quad
\mu=\left\{ \begin{array}{ll}
\beta \frac{\alpha}{\alpha-1} & \mbox{if $\alpha \not=1$}\\
\beta c_0 & \mbox{if $\alpha=1$}
\end{array} \right.$$
and $c_0= \int_0^{\infty}(\sin z-z 1_{\{z \leq 1\}})z^{-2}dz$. If $\alpha>1$, then $E(Y(B))=\mu|B|$.
\end{lemma}

\noindent {\bf Proof:} We first express the characteristic function (\ref{chf-Y-B}) of $Y(B)$ in Feller's canonical form (see Section XVII.2 of \cite{feller71}):
$$E(e^{iu Y(B)})= \exp \left\{ iub|B| +|B|\int_{\bR}\frac{e^{iuz}-1-iu \sin z}{z^2}M_{\alpha}(dz)\right\}$$
with $M_{\alpha}(dz)=z^2 \nu_{\alpha}(dz)$ and $b=\int_{\bR} (\sin z-z 1_{\{|z| \leq 1\}})\nu_{\alpha}(dz)$. Then the result follows from the calculations done in Example XVII.3.(g) of \cite{feller71}. $\Box$

\vspace{3mm}

From Lemma \ref{Y-indep-scattered} and Lemma \ref{Y(B)param}, it follows that
$$Z=\{Z(B)=Y(B)-\mu|B|; B \in \cB_{b}(\bR_{+} \times \bR^d)\}$$
is an $\alpha$-stable random measure, in the sense of Definition 3.3.1 of \cite{ST94}, with control measure $m(B)=\sigma^{\alpha}|B|$ and constant skewness intensity $\beta$. In particular, $Z(B)$ has a $S_{\alpha}(\sigma|B|^{1/\alpha},\beta,0)$ distribution.

We say that $Z$ is an {\bf $\alpha$-stable L\'evy noise}.
Coming back to the original construction (\ref{def-Y(B)}) of $Y(B)$ and noticing that
$$\mu|B|=-|B|\int_{\bR}z1_{\{|z| \leq 1\}} \nu_{\alpha}(dz)=-\sum_{j \geq 1}E(L_j(B)) \quad \mbox{if} \quad \alpha<1, \ \mbox{and}$$
$$\mu|B|=|B|\int_{\bR}z1_{\{|z| > 1\}} \nu_{\alpha}(dz)=E(L_0(B)) \quad \mbox{if} \quad \alpha>1,$$
it follows that $Z(B)$ can be represented as:
\begin{equation}
\label{def-Z-alpha-less}
Z(B)= \sum_{j \geq 0}L_j(B)=:\int_{B \times (\bR \verb2\2 \{0\})}z N(dt,dx,dz) \quad \mbox{if} \ \alpha<1,
\end{equation}
\begin{equation}
\label{def-Z-alpha-greater}
Z(B)=\sum_{j \geq 0}[L_j(B)-E(L_j(B))]=:\int_{B \times (\bR \verb2\2 \{0\})}z \widehat{N}(dt,dx,dz) \ \mbox{if} \ \alpha>1.
\end{equation}
Here $\widehat{N}$ is the compensated Poisson measure associated to $N$, i.e. $\widehat{N}(A)=N(A)-E(N(A))$ for any relatively compact set $A$ in $\bR_{+} \times \bR^d \times (\overline{\bR} \verb2\2 \{0\})$.

In the case $\alpha=1$, we will assume that $p=q$ so that $\nu_{\alpha}$ is symmetric around $0$, $E(L_j(B))=0$ for all $j \geq 1$, and $Z(B)$ admits the same representation as in the case $\alpha<1$.

\section{The linear equation}
\label{linear-eq-section}

As a preliminary investigation, we consider first equation (\ref{eq-noiseZ}) with $\sigma=1$:
\begin{equation}
\label{linear-eq}
Lu(t,x)=\dot{Z}(t,x), \quad t>0,x \in \cO
\end{equation}
with zero initial conditions and Dirichlet boundary conditions.  In this section $\cO$ is a bounded domain in $\bR^d$ or $\cO=\bR^d$. 

By definition, the process $\{u(t,x);t \geq 0,x \in \cO\}$ given by:
\begin{equation}
\label{def-sol-linear}
u(t,x)=\int_0^t \int_{\cO}G(t-s,x,y)Z(ds,dy)
\end{equation}
is a mild solution of (\ref{linear-eq}), provided that the stochastic integral on the right-hand side of (\ref{def-sol-linear}) is well-defined.


We define now the stochastic integral of a deterministic function $\varphi$: 
$$Z(\varphi)=\int_{0}^{\infty}\int_{\bR^d}\varphi(t,x)Z(dt,dx).$$

If $\varphi \in L^{\alpha}(\bR_{+} \times \bR^d)$, this can be defined by approximation with simple functions, as explained in Section 3.4 of \cite{ST94}.
The process $\{Z(\varphi); \varphi\in L^{\alpha}(\bR_{+} \times \bR^d)\}$ has jointly $\alpha$-stable finite dimensional distributions. In particular, each $Z(\varphi)$ has a $S_{\alpha}(\sigma_{\varphi},\beta,0)$-distribution with scale parameter: 
$$\sigma_{\varphi}=\sigma \left(\int_{0}^{\infty} \int_{\bR^d}|\varphi(t,x)|^{\alpha}dxdt \right)^{1/\alpha}.$$

More generally, a measurable function $\varphi: \bR_{+} \times \bR^d \to \bR$ is integrable with respect to $Z$ if there exists a sequence $(\varphi_n)_{n \geq 1}$ of simple functions such that $\varphi_n\to \varphi$ a.e., and for any $B \in \cB_{b}(\bR_{+} \times \bR^d)$, the sequence $\{Z(\varphi_n 1_{B})\}_{n}$ converges in probability (see \cite{RR89}).

The next results shows that condition $\varphi \in L^{\alpha}(\bR_{+} \times \bR^d)$ is also necessary for the integrability of $\varphi$ with respect to $Z$. Due to Lemma \ref{Y-indep-scattered}, this follows immediately from the general theory of stochastic integration with respect to independently scattered random measures developed in \cite{RR89}.

\begin{lemma}
\label{L-alpha-cond}
A deterministic function $\varphi$ is integrable with respect to $Z$ if and only if $\varphi \in L^{\alpha}(\bR_{+} \times \bR^d)$.
\end{lemma}

\noindent {\bf Proof:} We write the characteristic function of $Z(B)$ in the form used in \cite{RR89}:
$$E(e^{iuZ(B)})=\exp\left\{\int_{B}\left[iua+\int_{\bR}(e^{iuz}-1-iu \tau(z))\nu_{\alpha}(dz) \right]dtdx \right\}$$
with $a=\beta-\mu$, $\tau(z)=z$ if $|z| \leq 1$ and $\tau(z)={\rm sgn}(z)$ if $|z|>1$. By Theorem 2.7 of \cite{RR89}, $\varphi$ is integrable with respect to $Z$ if and only if
$$\int_{\bR_{+} \times \bR^d}|U(\varphi(t,x))|dtdx<\infty \quad \mbox{and} \quad\int_{\bR_{+} \times \bR^d}V(\varphi(t,x))dtdx<\infty$$
where
$U(y)=ay+\int_{\bR}(\tau(yz)-y\tau(z))\nu_{\alpha}(dz)$ and $V(y)=\int_{\bR}(1 \wedge |yz|^2)\nu_{\alpha}(dz)$. Direct calculations show that in our case, $U(y)=-\frac{\beta}{\alpha-1}y^{\alpha}$ if $\alpha\not=1$, $U(y)=0$ if $\alpha=1$, and $V(y)=\frac{2}{2-\alpha}y^{\alpha}$. $\Box$

\vspace{3mm}

The following result follows immediately from (\ref{def-sol-linear}) and Lemma \ref{L-alpha-cond}.

\begin{proposition}
Equation (\ref{linear-eq}) has a mild solution if and only if for any $t>0,x \in \cO$
\begin{equation}
\label{G-in-Lalpha}
I_{\alpha}(t)=\int_0^t \int_{\cO}G(s,x,y)^{\alpha}dyds<\infty.
\end{equation}
In this case, $\{u(t,x); t \geq 0,x\in \cO\}$ has jointly $\alpha$-stable finite-dimensional distributions. In particular, $u(t,x)$ has a $S_{\alpha}(\sigma I_{\alpha}(t)^{1/\alpha},\beta,0)$ distribution.
\end{proposition}



Condition (\ref{G-in-Lalpha}) can be easily verified in the case of several examples.

\begin{example} (Heat equation)
{\rm Let $L=\frac{\partial}{\partial t}-\frac{1}{2}\Delta$. Assume first that $\cO=\bR^d$. Then $G(t,x,y)=\overline{G}(t,x-y)$, where
\begin{equation}
\label{def-G-heat}
\overline{G}(t,x)=\frac{1}{(2\pi t)^{d/2}}\exp\left(-\frac{|x|^2}{2t}\right),
\end{equation}
and condition (\ref{G-in-Lalpha}) is equivalent to (\ref{heat-cond}). In this case,
$I_{\alpha}(t)=c_{\alpha,d}t^{d(1-\alpha)/2+1}$.
If $\cO$ is a bounded domain in $\bR^d$, then
$G(t,x,y) \leq \overline{G}(t,x-y)$ (see p. 74 of \cite{mueller-mytnik-stan06})
and condition (\ref{G-in-Lalpha}) is implied by (\ref{heat-cond}).
}
\end{example}

\begin{example} (Parabolic equations)
{\rm Let $L=\frac{\partial }{\partial t}-\cL$ where
\begin{equation}
\label{def-cL-operator}
\cL f(x)=\sum_{i,j=1}^{d}a_{ij}(x)\frac{\partial^2 f}{\partial x_i \partial x_j}(x)+\sum_{i=1}^{d}b_i(x)\frac{\partial f}{\partial x_i}(x)
\end{equation}
is the generator of a 
Markov process with values in $\bR^d$, without jumps (a diffusion).
Assume that $\cO$ is a bounded domain in $\bR^d$ or $\cO=\bR^d$. By Aronson estimate (see e.g. Theorem 2.6 of \cite{PZ07}), under some assumptions on the coefficients $a_{ij},b_i$, there exist some constants $c_1,c_2>0$ such that
\begin{equation}
\label{parabolic-bound}G(t,x,y) \leq c_1 t^{-d/2} \exp \left(-\frac{|x-y|^2}{c_2t}\right)
\end{equation}
for all $t>0$ and $x,y \in \cO$.
In this case, condition (\ref{G-in-Lalpha}) is implied by (\ref{heat-cond}).
}
\end{example}

\begin{example}(Heat equation with fractional power of the Laplacian)
{\rm Let $L=\frac{\partial}{\partial t}+(-\Delta)^{\gamma}$ for some $\gamma>0$. Assume that $\cO$ is a bounded domain in $\bR^d$ or $\cO=\bR^d$. Then
(see e.g. Appendix B.5 of \cite{PZ07})
\begin{equation}
\label{Green-function-fractional}
G(t,x,y)=\int_0^{\infty} \cG(s,x,y)g_{t,\gamma}(s)ds=\int_0^{\infty}\cG(t^{1/\gamma}s,x,y)g_{1,\gamma}(s)ds,
\end{equation}
where $\cG(t,x,y)$ is the fundamental solution of $\frac{\partial u}{\partial t}-\Delta u=0$ on $\cO$ and
$g_{t,\gamma}$ is the density of the measure $\mu_{t,\gamma}$, $(\mu_{t,\gamma})_{t \geq 0}$ being a convolution semigroup of measures on $[0,\infty)$
whose Laplace transform is given by:
$$\int_0^{\infty}e^{-us}g_{t,\gamma}(s)ds=\exp\left(-t u^{\gamma} \right), \quad \forall u>0.$$

Note that if $\gamma <1$, $g_{t,\gamma}$ is the density of $S_t$, where $(S_t)_{t \geq 0}$ is a $\gamma$-stable subordinator with L\'evy measure $\rho_{\gamma}(dx)=\frac{\gamma}{\Gamma(1-\gamma)}x^{-\gamma-1}1_{(0,\infty)}(x)dx$.

Assume first that $\cO=\bR^d$. Then $G(t,x,y)=\overline{G}(t,x-y)$, where
\begin{equation}
\label{Fourier-fractional}\overline{G}(t,x)=\int_{\bR^d} e^{i \xi \cdot x} e^{-t|\xi|^{2\gamma}}d\xi.
\end{equation}

If $\gamma<1$, then $\overline{G}(t,\cdot)$ is the density of $X_t$, $(X_t)_{t \geq 0}$ being a symmetric $(2\gamma)$-stable L\'evy process with values in $\bR^d$ defined by $X_t=W_{S_t}$, with $(W_t)_{t \geq 0}$ a Brownian motion in $\bR^d$ with variance $2$.  By Lemma \ref{fract-lemma} (Appendix B),
if $\alpha>1$, then (\ref{G-in-Lalpha}) holds if and only if
\begin{equation}
\label{heat-fractional-cond} \alpha<1+\frac{2\gamma}{d}.
\end{equation}

If $\cO$ is a bounded domain in $\bR^d$, then $G(t,x,y) \leq \overline{G}(t,x-y)$ (by Lemma 2.1 of \cite{mueller98}). In this case, if  $\alpha>1$, then (\ref{G-in-Lalpha}) is implied by (\ref{heat-fractional-cond}).
 }
\end{example}

\begin{example} (Cable equation in $\bR$)
{\rm Let $Lu=\frac{\partial u}{\partial t}-\frac{\partial^2 u}{\partial x^2}+u$ and $\cO= \bR$. Then $G(t,x,y)=\overline{G}(t,x-y)$, where
$$\overline{G}(t,x)=\frac{1}{\sqrt{4\pi t}} \exp\left(-\frac{|x|^2}{4t}-t\right),$$
and condition (\ref{G-in-Lalpha}) holds for any $\alpha \in (0,2)$.


}

\end{example}

\begin{example} (Wave equation in $\bR^d$ with $d=1,2$)
{\rm Let $L=\frac{\partial^2}{\partial t^2}-\Delta$ and $\cO=\bR^d$ with $d=1$ or $d=2$. Then $G(t,x,y)=\overline{G}(t,x-y)$, where
$$\overline{G}(t,x)=\frac{1}{2}1_{\{|x|<t\}} \quad \mbox{if} \ d=1$$
$$\overline{G}(t,x)=\frac{1}{2 \pi}\cdot \frac{1}{\sqrt{t^2-|x|^2}}1_{\{|x|<t\}} \quad \mbox{if} \ d=2.$$
Condition (\ref{G-in-Lalpha}) holds for any $\alpha \in (0,2)$.
In this case, $I_{\alpha}(t)=2^{-\alpha} t^2$ if $d=1$ and $I_{\alpha}(t)=\frac{(2\pi)^{1-\alpha}}{(2-\alpha)(3-\alpha)} t^{3-\alpha}$ if $d=2$.

}
\end{example}

\section{Stochastic integration}
\label{stoch-int-section}

In this section we construct a stochastic integral with respect $Z$ by generalizing the ideas of \cite{GM83} to the case of random fields. Unlike these authors, we do not assume that $Z(B)$ has a symmetric distribution, unless $\alpha =1$.

Let $\cF_t=\cF_t^N \vee \cN$ where $\cN$ is the $\sigma$-field of negligible sets in $(\Omega,\cF,P)$ and
$\cF_{t}^{N}$ is the $\sigma$-field generated by $N([0,s] \times A \times \Gamma)$ for all $s \in [0,t], A \in \cB_{b}(\bR^d)$ and for all Borel sets $\Gamma \subset \bR \verb2\2 \{0\}$ bounded away from $0$. Note that $\cF_{t}^{Z} \subset \cF_{t}^{N}$ where $\cF_{t}^Z$ is the $\sigma$-field generated by $Z([0,s] \times A),s \in [0,t],A \in \cB_{b}(\bR^d)$.

A process $X=\{X(t,x)\}_{t \geq 0,x \in \bR^d}$ is called {\em elementary} if it of the form
\begin{equation}
\label{elementary}
X(t,x)=1_{(a,b]}(t)1_{A}(x)Y
\end{equation}
where $0 \leq a<b$, $A \in \cB_{b}(\bR^d)$ and $Y$ is $\cF_{a}$-measurable and bounded. A {\em simple process} is a linear combination of elementary processes.
Note that any simple process $X$ can be written as:
\begin{equation}
\label{simple}
X(t,x)=1_{\{0\}}(t)Y_{0}(x)+\sum_{i=0}^{N-1}1_{(t_i,t_{i+1}]}(t)Y_{i}(x)
\end{equation}
with $0=t_0<t_1< \ldots<t_{N}<\infty$ and
$Y_i(x)=\sum_{j=1}^{m_i}1_{A_{ij}}(x)Y_{ij}$, where $(Y_{ij})_{j=1, \ldots,m_i}$ are $\cF_{t_i}$-measurable and $(A_{ij})_{j=1, \ldots,m_j}$ are disjoint sets in $\cB_{b}(\bR^d)$. Without loss of generality, we assume that $Y_0=0$.

We denote by $\cP$ the {\em predictable $\sigma$-field} on $\Omega \times \bR_{+} \times \bR^d$, i.e. the $\sigma$-field generated by all simple processes. We say that a process
$X=\{X(t,x)\}_{t \geq 0,x \in \bR^d}$ is {\em predictable} if the map $(\omega,t,x) \mapsto X(\omega,t,x)$ is $\cP$-measurable.

\begin{remark}
\label{leftcont-predictable}
{\rm One can show that the predictable $\sigma$-field $\cP$ is the $\sigma$-field generated by the class $\cC$ of processes $X$ such that $t \mapsto X(\omega,t,x)$ is left-continuous for any $\omega \in \Omega,x \in \bR^d$ and $(\omega,x) \mapsto X(\omega,t,x)$ is $\cF_t \times \cB(\bR^d)$-measurable for any $t>0$.
}

\end{remark}


Let $\cL_{\alpha}$ be the class of all predictable processes $X$ such that
$$\|X\|_{\alpha,T,B}^{\alpha}:=E\int_0^T \int_{B}|X(t,x)|^{\alpha}dxdt<\infty,$$
for all $T>0$ and $B \in \cB_b(\bR^d)$. Note that $\cL_{\alpha}$ is a linear space. 


Let $(E_k)_{k \geq 1}$ be an increasing sequence of sets in $\cB_{b}(\bR^d)$ such that $\bigcup_{k}E_k=\bR^d$. We define
$$\|X\|_{\alpha}=\sum_{k \geq 1}\frac{1 \wedge \|X\|_{\alpha,k,E_k}}{2^k} \quad \mbox{if} \quad \alpha>1,$$
$$\|X\|_{\alpha}^{\alpha}=\sum_{k \geq 1}\frac{1 \wedge \|X\|_{\alpha,k,E_k}^{\alpha}}{2^k} \quad \mbox{if} \quad \alpha \leq 1.$$

We identify two processes $X$ and $Y$ for which $\|X-Y\|_{\alpha}=0$, i.e. $X=Y$ $\nu$-a.e., where $\nu=Pdtdx$. In particular, we identify two processes $X$ and $Y$ if $X$ is a modification of $Y$, i.e. $X(t,x)=Y(t,x)$ a.s. for all $(t,x) \in \bR_{+} \times \bR^d$.

The space $\cL_{\alpha}$ becomes a metric space endowed with the metric $d_{\alpha}$:
$$d_{\alpha}(X,Y)=\|X-Y\|_{\alpha} \ \mbox{if} \ \alpha>1, \quad d_{\alpha}(X,Y)=\|X-Y\|_{\alpha}^{\alpha} \ \mbox{if} \ \alpha \leq 1.$$
 This follows using Minkowski's inequality if $\alpha > 1$, and the inequality $|a+b|^{\alpha} \leq |a|^{\alpha}+|b|^{\alpha}$ if $\alpha \leq 1$.

The following result can be proved similarly to Proposition 2.3 of \cite{walsh86}.

\begin{proposition}
\label{approx-simple}
For any $X \in \cL_{\alpha}$ there exists a sequence $(X_n)_{n \geq 1}$ of bounded simple processes such that $\|X_n-X\|_{\alpha} \to 0$ as $n \to \infty$.
\end{proposition}

By Proposition 5.7 of \cite{resnick07},  the $\alpha$-stable L\'evy process $\{Z(t,B)=Z([0,t] \times B);t \geq 0\}$ has a c\`adl\`ag modification, for any $B \in \cB_{b}(\bR^d)$. 
We work with these modifications. 
If $X$ is a simple process given by (\ref{simple}), we define
\begin{equation}
\label{int-simple-proc}
I(X)(t,B)=\sum_{i=0}^{N-1} \sum_{j=1}^{m_i}Y_{ij}Z( (t_i \wedge t,t_{i+1} \wedge t] \times (A_{ij} \cap B)).
\end{equation}
Note that for any $B \in \cB_{b}(\bR^d)$, $I(X)(t,B)$ is $\cF_t$-measurable for any $t \geq 0$, and $\{I(X)(t,B)\}_{t \geq 0}$ is c\`adl\`ag. We write
$$I(X)(t,B)=\int_0^t \int_{B} X(s,x)Z(ds,dx).$$

The following result will be used for the construction of the integral. This result generalizes Lemma 3.3 of \cite{GM83} to the case of random fields and non-symmetric measures $\nu_{\alpha}$.

\begin{theorem}
\label{max-ineq}
If $X$ is a bounded simple process then
\begin{equation}
\label{max-ineq-simple}
\sup_{\lambda>0}\lambda^{\alpha}P(\sup_{t \in [0,T]} |I(X)(t,B)|>\lambda) \leq c_{\alpha} E\int_0^T \int_{B}|X(t,x)|^{\alpha}dxdt,
\end{equation}
for any $T>0$ and $B \in \cB_{b}(\bR^d)$, where $c_{\alpha}$ is a constant depending only on $\alpha$.
\end{theorem}

\noindent {\bf Proof:} Suppose that $X$ is of the form (\ref{simple}). Since $\{I(X)(t,B)\}_{t \in [0,T]}$ is c\`adl\`ag, it is separable. Without loss of generality, we assume that its separating set $D$ can be written as $D=\cup_{n}F_n$ where $(F_n)_{n}$ is an increasing sequence of finite sets containing the points $(t_k)_{k=0, \ldots,N}$. Hence,
\begin{equation}
\label{max-ineq-step1}
P(\sup_{t \in[0,T]}|I(X)(t,B)|>\lambda)=\lim_{n \to \infty}P(\max_{t \in F_n}|I(X)(t,B)|>\lambda).
\end{equation}

Fix $n \geq 1$. Denote by $0=s_0<s_1< \ldots<s_m=T$ the points of the set $F_n$. Say $t_k=s_{i_k}$ for some $0=i_0<i_1<\ldots<i_N$.
Then each interval $(t_k,t_{k+1}]$ can be written as the union of some intervals of the form $(s_i,s_{i+1}]$:
\begin{equation}
\label{partition}
(t_k,t_{k+1}]=\bigcup_{i \in I_k}(s_i,s_{i+1}]
\end{equation}
where $I_k=\{i;i_k \leq i<i_{k+1}\}$. By (\ref{int-simple-proc}), for any $k=0, \ldots,N-1$ and $i \in I_k$,
$$I(X)(s_{i+1},B)-I(X)(s_i,B)=\sum_{j=1}^{m_k}Y_{kj} Z((s_i,s_{i+1}] \times (A_{kj} \cap B)).$$

For any $i \in I_k$, let $N_i=m_k$, and for any $j=1, \ldots, N_i$, define
$\beta_{ij}=Y_{kj}$, $H_{ij}=A_{kj}$ and $Z_{ij}=Z((s_i,s_{i+1}] \times (H_{ij}\cap B))$.
With this notation, we have:
$$I(X)(s_{i+1},B)-I(X)(s_i,B)=\sum_{j=1}^{N_i}\beta_{ij}Z_{ij} \quad \mbox{for all} \ i=0, \ldots,m.$$
Consequently, for any $l=1, \ldots,m$
\begin{equation}
\label{max-ineq-step2}I(X)(s_l,B)=\sum_{i=0}^{l-1}(I(X)(s_{i+1},B)-I(X)(s_i,B))=
\sum_{i=0}^{l-1}\sum_{j=1}^{N_i}\beta_{ij}Z_{ij}.
\end{equation}

Using (\ref{max-ineq-step1}) and  (\ref{max-ineq-step2}), it is enough to prove that for any $\lambda>0$,
\begin{equation}
\label{max-ineq-step3}P(\max_{l=0, \ldots,m-1} |\sum_{i=0}^{l}\sum_{j=1}^{N_i}\beta_{ij}Z_{ij}|>\lambda) \leq c_{\alpha} \lambda^{-\alpha}E \int_0^T \int_B |X(s,x)|^{\alpha}dx ds.
\end{equation}

First, note that
$$E \int_0^T \int_B |X(s,x)|^{\alpha}dx ds=\sum_{i=0}^{m-1}(s_{i+1}-s_i) \sum_{j=1}^{N_i}E|\beta_{ij}|^{\alpha} |H_{ij} \cap B|.$$
This follows from the definition (\ref{simple}) of $X$ and (\ref{partition}), since
$X(t,x)=\linebreak \sum_{i=0}^{N-1} \sum_{i \in I_k} 1_{(s_i,s_{i+1}]}(t) \sum_{j=1}^{N_i}\beta_{ij}1_{H_{ij}}(x)$.

We now prove (\ref{max-ineq-step3}). Let $W_i=\sum_{j=1}^{N_i}\beta_{ij}Z_{ij}$. For the event on the left-hand side, we consider its intersection with the event $\{\max_{0 \leq i \leq m-1}|W_i|>\lambda\}$ and its complement. Hence the probability of this event can be bounded by
$$\sum_{i=0}^{m-1}P(|W_i|>\lambda)+P(\max_{0 \leq l \leq m-1}|\sum_{i=0}^{l}W_i1_{\{|W_i| \leq \lambda\}}|>\lambda)=:I+II.$$
 We treat separately the two terms.

For the first term, we note that $\overline{\beta}_i=(\beta_{ij})_{1 \leq j \leq N_i}$ is $\cF_{s_i}$-measurable and $\overline{Z}_i=(Z_{ij})_{1 \leq j \leq N_i}$ is independent of $\cF_{s_i}$. By Fubini's theorem
$$I=\sum_{i=0}^{m-1}\int_{\bR^{N_i}} P(|\sum_{j=1}^{N_i}x_j Z_{ij}|>\lambda)P_{\overline{\beta}_i}(d\overline{x}),$$
where $\overline{x}=(x_j)_{1 \leq j \leq N_i}$ and $P_{\overline{\beta}_i}$ is the law of $\overline{\beta}_i$.

We examine the tail of $U_i=\sum_{j=1}^{N_i}x_j Z_{ij}$ for a fixed $\overline{x} \in \bR^{N_i}$. By Lemma \ref{Y(B)param}, $Z_{ij}$ has a $S_{\alpha}(\sigma (s_{i+1}-s_i)^{1/\alpha}|H_{ij} \cap B|^{1/\alpha},\beta,0)$ distribution. Since the sets $(H_{ij})_{1 \leq j \leq N_i}$ are disjoint, the variables $(Z_{ij})_{1 \leq j \leq N_i}$ are independent. Using elementary properties of the stable distribution (Properties 1.2.1 and 1.2.3 of \cite{ST94}), it follows that $U_i$ has a $S_{\alpha}(\sigma_i,\beta_i^*,0)$ distribution with parameters: $$\sigma_i^{\alpha} =\sigma^{\alpha}(s_{i+1}-s_i)\sum_{j=1}^{N_i}|x_j|^{\alpha}|H_{ij} \cap B|$$
$$\beta_i^*=\frac{\beta}{\sum_{j=1}^{N_i}|x_j|^{\alpha}|H_{ij} \cap B|} \sum_{j=1}^{N_i}{\rm sgn}(x_j)|x_j|^{\alpha}|H_{ij} \cap B|.$$

\noindent By Lemma \ref{lemmaA1} (Appendix A), there exists a constant $c_{\alpha}^*>0$ such that
\begin{equation}
\label{tail-Ui}
P(|U_i|>\lambda) \leq c_{\alpha}^{*}\lambda^{-\alpha}\sigma^{\alpha}(s_{i+1}-s_i)\sum_{j=1}^{N_i}|x_j|^{\alpha}|H_{ij} \cap B|
\end{equation}
for any $\lambda>0$. Hence
$$I \leq c_{\alpha}^*\lambda^{-\alpha}\sigma^{\alpha}\sum_{i=0}^{m-1}(s_{i+1}-s_i)
\sum_{j=1}^{N_i}E|\beta_{ij}|^{\alpha}|H_{ij} \cap B|=c_{\alpha}^*\lambda^{-\alpha}\sigma^{\alpha}E\int_0^T \int_B|X(s,x)|^{\alpha}dxds.$$

We now treat $II$. We consider three cases. For the first two cases we deviate from the original argument of \cite{GM83} since we do not require that $\beta=0$.

{\em Case 1.} $\alpha<1$.
Note that
\begin{equation}
\label{max-ineq-step4}II \leq P(\max_{0 \leq l \leq m-1}M_{l}>\lambda)
\end{equation}
where $\{M_l=\sum_{i=0}^{l}|W_i|1_{\{|W_i| \leq \lambda\}}, \cF_{s_{l+1}}; 0 \leq l \leq m-1\}$ is a submartingale. By the submartingale maximal inequality (Theorem 35.3 of \cite{billingsley95}),
\begin{equation}
\label{max-ineq-step5}P(\max_{0 \leq l \leq m-1}M_{l}>\lambda) \leq \frac{1}{\lambda}E(M_{m-1})=\frac{1}{\lambda}\sum_{i=0}^{m-1}E(|W_i|1_{|W_i| \leq \lambda}).
\end{equation}

Using the independence between $\overline{\beta}_i$ and $\overline{Z}_i$ it follows that
$$E[|W_i|1_{|W_i| \leq \lambda}]=\int_{\bR^{N_i}} E[|\sum_{j=1}^{N_i}x_jZ_{ij}|1_{\{|\sum_{j=1}^{N_i}x_jZ_{ij}| \leq \lambda \}}] P_{\overline{\beta}_i}(d\overline{x})$$
Let $U_i=\sum_{j=1}^{N_i}x_{j}Z_{ij}$. Using (\ref{tail-Ui}) and Remark \ref{Adam-remark} (Appendix \ref{appA-section}), we get:
$$E[|U_i|1_{\{|U_i| \leq \lambda\}}]
\leq c_{\alpha}^*\sigma^{\alpha}\frac{1}{1-\alpha}\lambda^{1-\alpha}(s_{i+1}-s_i)\sum_{j=1}^{N_i}|x_j|^{\alpha}|H_{ij} \cap B| .$$
Hence
\begin{equation}
\label{max-ineq-step6}E[|W_i|1_{|W_i| \leq \lambda}] \leq
c_{\alpha}^*\sigma^{\alpha}\frac{1}{1-\alpha}\lambda^{1-\alpha}
(s_{i+1}-s_i)\sum_{j=1}^{N_i}E|\beta_{ij}|^{\alpha}|H_{ij} \cap B|.
\end{equation}

From (\ref{max-ineq-step4}), (\ref{max-ineq-step5}) and (\ref{max-ineq-step6}), it follows that: $$II \leq c_{\alpha}^* \sigma^{\alpha} \frac{1}{1-\alpha}\lambda^{-\alpha}E\int_0^T \int_{B}|X(s,x)|^{\alpha}dxds.$$

{\em Case 2.} $\alpha>1$. We have
$$II \leq P(\max_{0 \leq l \leq m-1}|\sum_{i=0}^{l}X_i|>\lambda/2)+P(\max_{0 \leq l \leq m-1}Y_i>\lambda/2)=:II'+II'',$$
where $X_i=W_i 1_{\{|W_i| \leq \lambda\}}-E[W_i 1_{\{|W_i| \leq \lambda\}}|\cF_{s_i}]$ and
$Y_i=|E[W_i 1_{\{|W_i| \leq \lambda\}}|\cF_{s_i}]|$.

We first treat the term $II'$. Note that $\{M_l=\sum_{i=0}^{l}X_i,\cF_{s_{l+1}}; 0 \leq l \leq m-1\}$ is a zero-mean square integrable martingale, and
$$II'=P(\max_{0 \leq l \leq m-1}|M_l|>\lambda/2) \leq \frac{4}{\lambda^2} \sum_{i=0}^{m-1}E(X_i^2) \leq \frac{4}{\lambda^2}\sum_{i=0}^{m-1}E[W_i^2 1_{\{|W_i| \leq \lambda\}}].$$

\noindent 
Let $U_i=\sum_{j=1}^{N_i}x_j Z_{ij}$. Using (\ref{tail-Ui}) and Remark \ref{Adam-remark} (Appendix \ref{appA-section}), we get:
$$E[U_i^2 1_{\{|U_i| \leq \lambda\}}] 
\leq 2 c_{\alpha}^*\sigma^{\alpha}\frac{1}{2-\alpha}\lambda^{2-\alpha}(s_{i+1}-s_i)\sum_{j=1}^{N_i}|x_j|^{\alpha}|H_{ij} \cap B|.$$
As in Case 1, we obtain that:
\begin{equation}
\label{max-ineq-step8}E[W_i^2 1_{\{|W_i| \leq \lambda\}}] \leq c_{\alpha}^*\sigma^{\alpha} \frac{2}{2-\alpha}\lambda^{2-\alpha}(s_{i+1}-s_i)\sum_{j=1}^{N_i}E|\beta_{ij}|^{\alpha}|H_{ij} \cap B|,
\end{equation}
and hence
$$II' \leq 8 c_{\alpha}^* \sigma^{\alpha} \frac{1}{2-\alpha} \lambda^{-\alpha}E \int_0^T \int_{B}|X(s,x)|^{\alpha}dxds.$$

We now treat $II''$. Note that $\{N_l=\sum_{i=0}^{l}Y_i,\cF_{s_{l+1}}; 0 \leq l \leq m-1\}$ is a semimartingale and hence, by the submartingale inequality,
$$II'' \leq \frac{2}{\lambda}E(N_{m-1})=\frac{2}{\lambda}\sum_{i=0}^{m-1}E(Y_i).$$
To evaluate $E(Y_i)$, we note that for almost all $\omega \in \Omega$,
\begin{equation}
\label{max-ineq-step7}E[W_i 1_{\{|W_i| \leq \lambda\}}|\cF_{s_i}](\omega)=E[\sum_{j=1}^{N_i}\beta_{ij}(\omega)Z_{ij}
1_{\{|\sum_{j=1}^{N_i}\beta_{ij}(\omega)Z_{ij}| \leq \lambda\}}],
\end{equation}
due to the independence between $\overline{\beta}_i$ and $\overline{Z}_i$.
We let $U_i=\sum_{j=1}^{N_i}x_j Z_{ij}$ with $x_j=\beta_{ij}(\omega)$. Since $\alpha>1$, $E(U_i)=0$. Using (\ref{tail-Ui}) and Remark \ref{Adam-remark}, we obtain
\begin{eqnarray*}
|E[U_i 1_{\{|U_i| \leq \lambda\}}]|&=&|E[U_i 1_{\{|U_i|>\lambda\}}]| \leq E[|U_i| 1_{\{|U_i| > \lambda\}}] \\
&\leq & c_{\alpha}^* \sigma^{\alpha}\frac{\alpha}{\alpha-1}\lambda^{1-\alpha} (s_{i+1}-s_i) \sum_{j=1}^{N_i}|x_j|^{\alpha}|H_{ij} \cap B|.
\end{eqnarray*}
Hence,
$E(Y_i) \leq c_{\alpha}^* \sigma^{\alpha}\frac{\alpha}{\alpha-1}\lambda^{1-\alpha} (s_{i+1}-s_i) \sum_{j=1}^{N_i}E|\beta_{ij}|^{\alpha}|H_{ij} \cap B|$ and
$$II'' \leq c_{\alpha}^* \sigma^{\alpha}\frac{2\alpha}{\alpha-1}\lambda^{-\alpha}E\int_0^T \int_{B}|X(t,x)|^{\alpha}dxdt.$$

{\em Case 3.} $\alpha=1$. In this case we assume that $\beta=0$. Hence $U_i=\sum_{j=1}^{N_i}x_j Z_{ij}$ has a symmetric distribution for any $\overline{x} \in \bR^{N_i}$. Using (\ref{max-ineq-step7}), it follows that $E[W_i 1_{\{|W_i| \leq \lambda\}}|{\cF}_{s_{i}}]=0$ a.s. for all $i=0, \ldots,m-1$. Hence, $\{M_l=\sum_{i=0}^{l}W_i 1_{\{|W_i| \leq \lambda\}}, \cF_{s_{l+1}}; 0 \leq l \leq m-1\}$ is a zero-mean square integrable martingale. By the martingale maximal inequality,
$$II \leq \frac{1}{\lambda^2}E[M_{m-1}^2]=\frac{1}{\lambda^2}\sum_{i=0}^{m-1}E[W_i^2 1_{\{|W_i| \leq \lambda\}}].$$
The result follows using (\ref{max-ineq-step8}).
$\Box$

\vspace{3mm}

We now proceed to the construction of the stochastic integral. If $Y=\{Y(t)\}_{t \geq 0}$ is a jointly measurable random process, we define:
$$\|Y\|_{\alpha,T}^{\alpha}=\sup_{\lambda>0}\lambda^{\alpha}P(\sup_{t \in [0,T]}|Y(t)|>\lambda).$$

Let $X \in \cL_{\alpha}$ be arbitrary. By Proposition \ref{approx-simple}, there exists a sequence $(X_n)_{n \geq 1}$ of simple functions such that $\|X_n-X\|_{\alpha} \to 0$ as $n \to \infty$. Let $T>0$ and $B \in \cB_{b}(\bR^d)$ be fixed. By linearity of the integral and Theorem \ref{max-ineq},
\begin{equation}
\label{Cauchy-seq}\|I(X_n)(\cdot,B)-I(X_m)(\cdot,B)\|_{\alpha,T}^{\alpha} \leq c_{\alpha}\|X_n-X_m\|_{\alpha,T,B}^{\alpha} \to 0
\end{equation}
as $n,m \to \infty$. In particular, the sequence $\{I(X_n)(\cdot,B)\}_{n}$ is Cauchy in probability in the space $D[0,T]$ equipped with the sup-norm. Therefore, there exists a random element $Y(\cdot,B)$ in $D[0,T]$ such that for any $\lambda>0$,
$$P(\sup_{t \in [0,T]}|I(X_n)(t,B)-Y(t,B)|>\lambda) \to 0.$$
Moreover, there exists a subsequence $(n_k)_k$ such that
$$\sup_{t \in [0,T]}|I(X_{n_k})(t,B)-Y(t,B)| \to 0 \quad {\rm a.s.}$$
as $k \to \infty$. Hence $Y(t,B)$ is $\cF_{t}$-measurable for any $t \in [0,T]$. The process $Y(\cdot,B)$ does not depend on the sequence $(X_n)_n$ and can be extended to a c\`adl\`ag process on $[0,\infty)$, which is unique up to indistinguishability. We denote this extension by $I(X)(\cdot,B)$ and we write
$$I(X)(t,B)=\int_{0}^{t}\int_{B}X(s,x)Z(ds,dx).$$
If $A$ and $B$ are disjoint sets in $\cB_{b}(\bR^d)$, then
\begin{equation}
\label{additivity}I(X)(t,A \cup B)=I(X)(t,A)+I(X)(t,B) \quad \mbox{a.s.}
\end{equation}

\begin{lemma}
\label{max-ineq-L}
Inequality (\ref{max-ineq-simple}) holds for any $X \in \cL_{\alpha}$.
\end{lemma}

\noindent {\bf Proof:} Let $(X_n)_n$ be a sequence of simple functions such that $\|X_n-X\|_{\alpha} \to 0$. For fixed $B$, we denote $I(X)=I(X)(\cdot,B)$. We let $\|\cdot\|_{\infty}$ be the sup norm on $D[0,T]$.
For any $\varepsilon>0$, we have:
$$P(\|I(X)\|_{\infty}>\lambda) \leq P(\|I(X)-I(X_n)\|_{\infty}>\lambda \varepsilon) +P(\|I(X_n)\|_{\infty}>\lambda (1-\varepsilon)).$$
Multiplying by $\lambda^{\alpha}$, and using Theorem \ref{max-ineq}, we obtain:
$$\sup_{\lambda>0}\lambda^{\alpha}P(\|I(X)\|_{\infty}>\lambda) \leq \varepsilon^{-\alpha}
\sup_{\lambda>0}\lambda^{\alpha}P(\|I(X)-I(X_n)\|_{\infty}>\lambda)+
(1-\varepsilon)^{-\alpha}c_{\alpha}\|X_n\|_{\alpha,T,B}^{\alpha}.$$
Let $n \to \infty$. Using (\ref{Cauchy-seq}) one can prove that
$\sup_{\lambda>0}\lambda^{\alpha}P(\|I(X_n)-I(X)\|_{\infty}>\lambda) \to 0$. We obtain that $\sup_{\lambda>0}\lambda^{\alpha}P(\|I(X)\|_{\infty}>\lambda) \leq (1-\varepsilon)^{-\alpha}c_{\alpha}\|X\|_{\alpha,T,B}^{\alpha}$. The conclusion follows
letting $\varepsilon \to 0$. $\Box$

\vspace{3mm}

For an arbitrary Borel set $\cO \subset \bR^d$ (possibly $\cO=\bR^d$), we assume in addition, that $X \in \cL_{\alpha}$ satisfies the condition:
\begin{equation}
\label{finite-int-Rd}
E\int_0^T \int_{\cO}|X(t,x)|^{\alpha}dxdt<\infty \quad \mbox{for all} \ T>0.
\end{equation}
Then we can define $I(X)(\cdot,\cO)$ as follows. Let $\cO_k=\cO \cap E_k$ where $(E_k)_k$ is an increasing sequence of sets in $\cB_{b}(\bR^d)$ such that $\bigcup_{k}E_k=\bR^d$. By (\ref{additivity}), Lemma
\ref{max-ineq-L} and (\ref{finite-int-Rd}),
$$\sup_{\lambda>0}\lambda^{\alpha}P(\sup_{t \leq T}|I(X)(t,\cO_k)-I(X)(t,\cO_l)|>\lambda) \leq c_{\alpha}E\int_0^T \int_{\cO_k\verb2\2 \cO_l}|X(t,x)|^{\alpha}dxdt \to 0$$
as $k,l \to \infty$. This shows that $\{I(X)(\cdot,\cO_k)\}_{k}$ is a Cauchy sequence in probability in the space $D[0,T]$ equipped with the sup norm. We denote by $I(X)(\cdot,\cO)$ its limit. As above, this process can be extended to $[0,\infty)$ and $I(X)(t,\cO)$ is $\cF_t$-measurable for any $t>0$. We denote
$$I(X)(t,\cO)=\int_0^t \int_{\cO}X(s,x)Z(ds,dx).$$
Similarly to Lemma \ref{max-ineq-L}, one can prove that for any $X \in \cL_{\alpha}$ satisfying (\ref{finite-int-Rd}),
$$\sup_{\lambda>0}\lambda^{\alpha}P(\sup_{t \leq T}|I(X)(t,\cO)|>\lambda) \leq c_{\alpha}E\int_0^T \int_{\cO}|X(t,x)|^{\alpha}dxdt.$$

\section{The truncated noise}
\label{truncated-noise-section}

For the study of non-linear equations, we need to develop a theory of stochastic integration with respect to another process $Z_K$ which is defined by removing from $Z$ the jumps whose modulus exceed a fixed value $K>0$. More precisely, for any $B \in \cB_{b}(\bR_{+} \times \bR^d)$, we define
\begin{equation}
\label{def-ZK-alpha-less}
Z_{K}(B)=
\int_{B \times \{0<|z| \leq K\}}zN(ds,dx,dz) \quad \mbox{if} \quad \alpha \leq 1
\end{equation}
\begin{equation}
\label{def-ZK-alpha-greater}Z_{K}(B)=
\int_{B \times \{0<|z| \leq K\}}z\widehat{N}(ds,dx,dz) \quad \mbox{if} \quad \alpha > 1.
\end{equation}

We treat separately the cases $\alpha \leq 1$ and $\alpha>1$.

\subsection{The case $\alpha \leq 1$}

Note that $\{Z_{K}(B); B \in \cB_{b}(\bR_{+} \times \bR^d)\}$ is an independently scattered random measure on $\bR_{+} \times \bR^d$ with characteristic function given by:
$$E(e^{iu Z_K(B)})=\exp \left\{|B|\int_{|z| \leq K}(e^{iuz}-1)\nu_{\alpha}(dz) \right\} \quad \forall \ u \in \bR.$$

We first examine the tail of $Z_{K}(B)$.

\begin{lemma}
\label{tail-ZKB-lemma}
For any set $B \in \cB_{b}(\bR_{+} \times \bR^d)$,
\begin{equation}
\label{tail-ZKB}
\sup_{\lambda>0}\lambda^{\alpha}P(|Z_K(B)|>\lambda) \leq r_{\alpha}|B|
\end{equation}
where $r_{\alpha}>0$ a constant depending only on $\alpha$ (given by Lemma \ref{lemmaA2}).
\end{lemma}

\noindent {\bf Proof:}
This follows from Example 3.7 of \cite{GM83}. We denote by $\nu_{\alpha,K}$ the restriction of $\nu_{\alpha}$ to $\{z \in \bR;0<|z| \leq K\}$. Note that
$$\nu_{\alpha,K}(\{z \in \bR;|z|>t\})=
\left\{
\begin{array}{ll}
t^{-\alpha}-K^{-\alpha} & \mbox{if $0<t \leq K$} \\
0 & \mbox{if $t>K$}
\end{array}
\right.$$
and hence $\sup_{t>0}t^{\alpha}\nu_{\alpha,K}(\{z \in \bR; |z|>t\})=1$.
Next we observe that we do not need to assume that the measure $\nu_{\alpha,K}$ is symmetric since we use a modified version of Lemma 2.1 of \cite{GM82} given by Lemma \ref{lemmaA2} (Appendix A). $\Box$

\vspace{3mm}

In fact, since the tail of $\nu_{\alpha,K}$ vanishes if $t>K$, we can obtain another estimate for the tail of $Z_{K}(B)$ which, together with (\ref{tail-ZKB}), will allow us to control its $p$-th moment for $p \in (\alpha,1)$. This new estimate is given below.

\begin{lemma}
\label{tail-ZKB-power1} If $\alpha<1$, then
$$P(|Z_K(B)|>u) \leq \frac{\alpha}{1-\alpha}K^{1-\alpha}|B|u^{-1} \quad \mbox{for all} \ u>K.$$
If $\alpha=1$, then
$P(|Z_K(B)|>u) \leq K|B|u^{-2}$ for all $u>K$.
\end{lemma}

\noindent {\bf Proof:} We use the same idea as in Example 3.7 of \cite{GM83}. For each $k \geq 1$, let $Z_{k,K}(B)$ be a random variable with characteristic function:
$$E(e^{iuZ_{k,K}(B)})=\exp\left\{|B|\int_{\{k^{-1}<|z| \leq K\}} (e^{iuz}-1)\nu_{\alpha}(dz) \right\}.$$
Since $\{Z_{k,K}(B)\}_{k}$ converges in distribution to $Z_{K}(B)$, it suffices to prove the lemma for $Z_{k,K}(B)$. Let $\mu_{k}$ be the restriction of $\nu_{\alpha}$ to $\{z; k^{-1}<|z|\leq K\}$. Since $\mu_{k}$ is finite, $Z_{k,K}(B)$ has a compound Poisson distribution with
\begin{equation}
\label{cPoisson-tail}P(|Z_{k,K}(B)|>u)=e^{-|B|\mu_{k}(\bR)}\sum_{n \geq 0}\frac{|B|^n}{n!} \mu_{k}^{*n}(\{z;|z|>u\}).
\end{equation}
where $\mu_{k}^{*n}$ denotes the $n$-fold convolution. Note that
$$\mu_{k}^{*n}(\{z;|z|>u\})=[\mu_k(\bR)]^nP(|\sum_{i=1}^{n}\eta_i|>u),$$
where $(\eta_i)_{i \geq 1}$ are i.i.d. random variables with law $\mu_k/\mu_k(\bR)$.

Assume first that $\alpha<1$. To compute $P(|\sum_{i=1}^{n}\eta_i|>u)$ we consider the intersection with the event $\{\max_{1 \leq i \leq n}|\eta_i|>u\}$ and its complement. Note that
$P(|\eta_i|>u)=0$ for any $u>K$. Using this fact and Markov's inequality, we obtain that for any $u>K$,
$$P(|\sum_{i=1}^{n}\eta_i|>u) \leq P(|\sum_{i=1}^{n}\eta_i 1_{\{|\eta_i| \leq u\}}|>u) \leq \frac{1}{u}\sum_{i=1}^{n}E(|\eta_i|1_{\{\{|\eta_i| \leq u \}}).$$

\noindent Note that $P(|\eta_i|>s) \leq (s^{-\alpha}-K^{-\alpha})/\mu_{k}(\bR)$ if $s \leq K$. Hence, for any $u>K$
$$E(|\eta_{i}| 1_{\{|\eta_i| \leq u\}}) \leq \int_0^u P(|\eta_i|>s)ds=\int_0^K P(|\eta_i|>s)ds \leq \frac{1}{\mu_{k}(\bR)}\frac{\alpha}{1-\alpha}K^{1-\alpha}.$$
Combining all these facts, we get: for any $u>K$
$$\mu_{k}^{*n}(\{z;|z|>u\}) \leq [\mu_k(\bR)]^{n-1}  \frac{\alpha}{1-\alpha}K^{1-\alpha}n u^{-1},$$
and the conclusion follows from (\ref{cPoisson-tail}).

Assume now that $\alpha=1$. In this case, $E(\eta_i 1_{\{|\eta_i| \leq u\}})=0$ since $\eta_i$ has a symmetric distribution. Using Chebyshev's inequality this time, we obtain:
$$P(|\sum_{i=1}^{n}\eta_i|>u) \leq P(|\sum_{i=1}^{n}\eta_i 1_{\{|\eta_i| \leq u\}}|>u) \leq \frac{1}{u^2}\sum_{i=1}^{n}E(\eta_i^2 1_{\{\{|\eta_i| \leq u \}}).$$
The result follows as above using the fact that for any $u>K$,
$$E(\eta_i^2 1_{\{|\eta_i| \leq u\}}) \leq 2 \int_0^u s P(|\eta_i|>s)ds=2\int_0^K sP(|\eta_i|>s)ds \leq \frac{1}{\mu_k(\bR)}K.$$
$\Box$

\begin{lemma}
\label{momentZK-alpha-less1}
If $\alpha<1$ then $$E|Z_K(B)|^p \leq C_{\alpha,p}K^{p-\alpha}|B| \quad \mbox{for any} \ p \in (\alpha,1),$$ where $C_{\alpha,p}$ is a constant depending on $\alpha$ and $p$. If $\alpha=1$, then $$E|Z_K(B)|^p \leq C_{p}K^{p-1}|B| \quad \mbox{for any} \ p \in (1,2),$$ where $C_{p}$ is a constant depending on $p$.
\end{lemma}

\noindent {\bf Proof:} Note that
$$E|Z_K(B)|^{p}=\int_0^{\infty}P(|Z_K(B)|^p>t)dt=p\int_0^{\infty}P(|Z_K(B)|>u)u^{p-1}du.$$
We consider separately the integrals for $u \leq K$ and $u>K$. For the first integral we use (\ref{tail-ZKB}):
$$\int_0^K P(|Z_K(B)|>u)u^{p-1}du \leq r_{\alpha}|B|\int_0^K u^{-\alpha+p-1}du=r_{\alpha}|B|\frac{1}{p-\alpha}K^{p-\alpha}.$$
For the second one we use Lemma \ref{tail-ZKB-power1}: if $\alpha<1$ then
$$\int_{K}^{\infty}P(|Z_K(B)|>u)u^{p-1}du \leq  \frac{\alpha}{1-\alpha}K^{1-\alpha}|B|\int_{K}^{\infty}u^{p-2}du
= \frac{\alpha}{(1-\alpha)(1-p)}|B|K^{p-\alpha},$$
and if $\alpha=1$, then
$$\int_{K}^{\infty}P(|Z_K(B)|>u)u^{p-1}du \leq  K|B| \int_{K}^{\infty}u^{p-3}du=|B|\frac{1}{2-p}K^{p-1}.$$
$\Box$

\vspace{3mm}

We now proceed to the construction of the stochastic integral with respect to $Z_K$. For this, we use the same method as for $Z$. Note that $\cF_t^{Z_{K}} \subset \cF_{t}$, where $\cF_t^{Z_K}$ is the $\sigma$-field generated by $Z_{K}([0,s] \times A)$ for all $s \in [0,t]$ and $A \in \cB_{b}(\bR^d)$. For any $B \in \cB_{b}(\bR^d)$, we will work with a c\`adl\`ag modification of the L\'evy process $\{Z_{K}(t,B)=Z_K([0,t] \times B);t \geq 0\}$.

If $X$ is a simple process given by (\ref{simple}), we define
$$I_K(X)(t,B)=\int_0^t \int_{B}X(s,x)Z_{K}(ds,dx)$$
by the same formula (\ref{int-simple-proc}) with $Z$ replaced by $Z_K$. The following result shows that $I_K(X)(t,B)$ has the same tail behavior as $I(X)(t,B)$.

\begin{proposition}
\label{max-ineq-alpha-less-1}
If $X$ is a bounded simple process then
\begin{equation}
\label{max-ineqZK}
\sup_{\lambda>0}\lambda^{\alpha}P(\sup_{t \in [0,T]} |I_K(X)(t,B)|>\lambda) \leq d_{\alpha} E\int_0^T \int_{B}|X(t,x)|^{\alpha}dxdt,
\end{equation}
for any $T>0$ and $B \in \cB_{b}(\bR^d)$, where $d_{\alpha}$ is a constant depending only on $\alpha$.
\end{proposition}

\noindent {\bf Proof:} As in the proof of Theorem \ref{max-ineq}, it is enough to prove that
$$P(\max_{l=0, \ldots,m-1} |\sum_{i=0}^{l}\sum_{j=1}^{N_i}\beta_{ij}Z_{ij}^*|>\lambda) \leq d_{\alpha} \lambda^{-\alpha} \sum_{i=0}^{m-1}(s_{i+1}-s_i) \sum_{j=1}^{N_i}E|\beta_{ij}|^{\alpha}|H_{ij} \cap B|,$$
where $Z_{ij}^{*}=Z_{K}((s_{i},s_{i+1}] \times (H_{ij} \cap B))$. This reduces to showing that $U_i^*=\sum_{j=1}^{N_i}x_jZ_{ij}^*$ satisfies an inequality similar to (\ref{tail-Ui}) for any $\overline{x} \in \bR^{N_i}$, i.e.
\begin{equation}
\label{max-ineqZK-step1}
P(|U_{i}^*|>\lambda) \leq d_{\alpha}^* \lambda^{-\alpha}(s_{i+1}-s_i) \sum_{j=1}^{N_i}|x_j|^{\alpha}|H_{ij} \cap B|,
\end{equation}
for any $\lambda>0$, for some $d_{\alpha}^*>0$. We first examine the tail of $Z_{ij}^*$. By (\ref{tail-ZKB}),
$$P(|Z_{ij}^*|>\lambda) \leq r_{\alpha}(s_{i+1}-s_i)K_{ij}\lambda^{-\alpha}.$$
where $K_{ij}=|H_{ij} \cap B|$. Letting $\eta_{ij}=K_{ij}^{-1/\alpha}Z_{ij}^*$, we obtain that for any $u>0$,
$$P(|\eta_{ij}|>u) \leq r_{\alpha}(s_{i+1}-s_i) u^{-\alpha} \quad \forall j=1,\ldots,N_i.$$
By Lemma \ref{lemmaA2} (Appendix A), it follows that for any $\lambda>0$,
$$P(|\sum_{j=1}^{N_i}b_j \eta_{ij}|>\lambda) \leq r_{\alpha}^2(s_{i+1}-s_i)\sum_{j=1}^{N_i}|b_j|^{\alpha}\lambda^{-\alpha},$$
for any sequence $(b_j)_{j=1, \ldots,N_i}$ of real numbers. Inequality (\ref{max-ineqZK-step1}) (with $d_{\alpha}^*=r_{\alpha}^2$) follows by applying this to $b_j=x_jK_{ij}^{1/\alpha}$. $\Box$

\vspace{3mm}

In view of the previous result and Proposition \ref{approx-simple}, for any process $X \in \cL_{\alpha}$ we can construct the integral $$I_K(X)(t,B)=\int_0^t \int_{B}X(s,x)Z_{K}(ds,dx)$$ in the same manner as $I(X)(t,B)$, and this integral satisfies (\ref{max-ineqZK}). If in addition the process $X \in \cL_{\alpha}$ satisfies (\ref{finite-int-Rd}), then we can define the integral $I_K(X)(t,\cO)$ for an arbitrary Borel set $\cO \subset \bR^d$ (possibly $\cO=\bR^d$). This integral will satisfy an inequality similar to (\ref{max-ineqZK}) with $B$ replaced by $\cO$.

The appealing feature of $I_{K}(X)(t,B)$ is that we can control its moments, as shown by the next result.

\begin{theorem}
\label{moment-ineq}
If $\alpha<1$, then for any $p \in (\alpha,1)$ and for any $X \in \cL_{p}$,
\begin{equation}
\label{p-moment-ineq}E|I_K(X)(t,B)|^p \leq C_{\alpha,p}K^{p-\alpha}E\int_0^t \int_{B}|X(s,x)|^p dxds,
\end{equation}
for any $t>0$ and $B \in \cB_{b}(\bR^d)$, where $C_{\alpha,p}$ is a constant depending on $\alpha,p$. If $\cO \subset \bR^d$ is an arbitrary Borel set, and we assume in addition, that the process $X\in \cL_p$ satisfies:
\begin{equation}
\label{p-finite-int-Rd}E\int_0^T \int_{\cO}|X(s,x)|^pdxds<\infty \quad \forall \ T>0,
\end{equation}
then inequality (\ref{p-moment-ineq}) holds with $B$ replaced by $\cO$.
\end{theorem}

\noindent {\bf Proof:} {\em Step 1.} Suppose that $X$ is an elementary process of the form (\ref{elementary}). Then
$I_K(X)(t,B)=YZ_{K}(H)$ where $H=(t \wedge a,t \wedge b] \times (A \cap B)$. Note that $Z_{K}(H)$ is independent of $\cF_{a}$. Hence, $Z_{K}(H)$ is independent of $Y$. Let $P_{Y}$ denote the law of $Y$. By Fubini's theorem,
\begin{eqnarray*}
E|YZ_{K}(H)|^p &=& p \int_0^{\infty}P(|YZ_{K}(H)|>u)u^{p-1}du \\
&=& p \int_{\bR}  \left(\int_0^{\infty}P(|yZ_{K}(H)|>u) u^{p-1} du \right)P_{Y}(dy).
\end{eqnarray*}
We evaluate the inner integral. We split this integral into two parts, for $u \leq K|y|$, respectively $u>K|y|$. For the first integral, we use (\ref{tail-ZKB}). For the second one, we use Lemma \ref{tail-ZKB-power1}. Therefore, the inner integral is bounded by:
$$r_{\alpha}|y|^{\alpha}|H| \int_{0}^{K|y|}u^{-\alpha+p-1}du+ \frac{\alpha}{1-\alpha} |y| K^{1-\alpha}|H|\int_{K|y|}^{\infty}u^{p-2}du=C_{\alpha,p}'K^{p-\alpha}|y|^p |H|$$
and
\begin{eqnarray*}
E|YZ_{K}(H)|^p  & \leq & 
pC_{\alpha,p}'K^{p-\alpha}|H|E|Y|^p
= C_{\alpha,p} K^{p-\alpha}E\int_0^t \int_{B}|X(s,x)|^{p} dx ds.
\end{eqnarray*}

{\em Step 2.} Suppose now that $X$ is a simple process of the form (\ref{simple}). Then $X(t,x)=\sum_{i=0}^{N-1}\sum_{j=1}^{m_i}X_{ij}(t,x)$ where $X_{ij}(t,x)=1_{(t_i,t_{i+1}]}(t)1_{A_{ij}}(x)Y_{ij}$.

Using the linearity of the integral, the inequality $|a+b|^p \leq |a|^p+|b|^p$, and the result obtained in Step 1 for the elementary processes $X_{ij}$, we get:
\begin{eqnarray*}
& & E|I_K(X)(t,B)|^p \leq E \sum_{i=0}^{N-1}\sum_{j=1}^{m_i}|I_{K}(X_{ij})(t,B)|^p \leq \\
& &  C_{\alpha,p}K^{p-\alpha} E \sum_{i=0}^{N-1}\sum_{j=1}^{m_i}\int_0^t \int_{B}|X_{ij}(s,x)|^pdxds =C_{\alpha,p}K^{p-\alpha}E \int_0^t \int_{B}|X(s,x)|^pdxds.
\end{eqnarray*}

{\em Step 3.} Let $X \in \cL_{p}$ be arbitrary. By Proposition \ref{approx-simple}, there exists a sequence $(X_n)_n$ of bounded simple processes such that $\|X_n-X\|_{p} \to 0$. Since $\alpha<p$, it follows that $\|X_n-X\|_{\alpha} \to 0$. By the definition of $I_K(X)(t,B)$ there exists a subsequence $\{n_k\}_k$ such that $\{I_{K}(X_{n_k})(t,B)\}_k$ converges to $I_K(X)(t,B)$ a.s. Using Fatou's lemma and the result obtained in Step 2 (for the simple processes $X_{n_k}$), we get:
\begin{eqnarray*}
E|I_K(X)(t,B)|^p & \leq & \liminf_{k \to \infty} E|I_{K}(X_{n_k})(t,B)|^p \\
& \leq & C_{\alpha,p}K^{p-\alpha} \liminf_{k \to \infty} E \int_0^t \int_B |X_{n_k}(s,x)|^p dxds \\
&=& C_{\alpha,p}K^{p-\alpha} E \int_0^t \int_B |X(s,x)|^p dxds.
\end{eqnarray*}

{\em Step 4.} Suppose that $X \in  \cL_p$ satisfies (\ref{p-finite-int-Rd}). Let $\cO_{k}=\cO \cap E_k$ where $(E_k)_k$ is an increasing sequence of sets in $\cB_b(\bR^d)$ such that $\bigcup_{k \geq 1}E_k=\bR^d$. By the definition of $I_K(X)(t,\cO)$, there exists a subsequence $(k_i)_{i}$ such that $\{I_{K}(X)(t,\cO_{k_i})\}_{i}$ converges to $I_K(X)(t,\cO)$ a.s. Using Fatou's lemma, the result obtained in Step 3 (for $B=\cO_{k_i}$) and the monotone convergence theorem, we get:
\begin{eqnarray*}
E|I_K(X)(t,\cO)|^p & \leq & \liminf_{i \to \infty} E|I_{K}(X)(t,\cO_{k_i})|^p \\
& \leq & C_{\alpha,p}K^{p-\alpha} \liminf_{i \to \infty} E \int_0^t \int_{\cO_{k_i}} |X(s,x)|^p dxds \\
&=& C_{\alpha,p}K^{p-\alpha} E \int_0^t \int_{\cO} |X(s,x)|^p dxds.
\end{eqnarray*}
$\Box$

\begin{remark}
{\rm Finding a similar moment inequality for the case $\alpha=1$ and $p \in (1,2)$ remains an open problem. The argument used in Step 2 above relies on the fact that $p<1$. Unfortunately, we could not find another argument to cover the case $p>1$.}
\end{remark}

\subsection{The case $\alpha>1$}
\label{tr-noise-alpha-greater-section}

In this case, the construction of the integral with respect to $Z_K$ relies on an integral with respect to $\widehat{N}$ which exists in the literature.
We recall briefly the definition of this integral. For more details, see Section 1.2.2 of \cite{bie98}, Section 24.2 of \cite{truman-wu06} or Section 8.7 of \cite{PZ07}.


Let $\bE=\bR^d \times (\bR \verb2\2\{0\})$ endowed with the measure $\mu(dx,dz)=dx \nu_{\alpha}(dz)$ and $\cB_{b}(\bE)$ be the class of bounded Borel sets in $\bE$. For a simple process $Y=\{Y(t,x,z);t \geq 0, (x,z) \in \bE\}$, the integral $I^{\widehat{N}}(Y)(t,B)$  is defined in the usual way, for any $t>0,B \in \cB_{b}(\bE)$. The process $I^{\widehat{N}}(Y)(\cdot,B)$ is a (c\`adl\`ag) zero-mean square-integrable martingale with quadratic variation
$$[I^{\widehat{N}}(Y)(\cdot,B)]_t=\int_0^t \int_{B}|Y(s,x,z)|^2 N(ds,dx,dz)$$
and predictable quadratic variation
$$\langle I^{\widehat{N}}(Y)(\cdot,B) \rangle_t=\int_0^t \int_{B}|Y(s,x,z)|^2 \nu_{\alpha}(dz)dx ds.$$
By approximation, this integral can be extended to the class of all $\cP \times \cB(\bR \verb2\2 \{0\})$-measurable processes $Y$ such that, for any $T>0$ and $B \in \cB_{b}(\bE)$
$$\|Y \|_{2,T,B}^2:=E \int_0^T \int_{B}|Y(s,x,z)|^2 \nu_{\alpha}(dz) dx ds <\infty.$$
The integral is a martingale with the same quadratic variations as above, and has the isometry property:
$E|I^{\widehat{N}}(Y)(t,B)|^2=\|Y\|_{2,T,B}^2$.
If in addition, $\|Y\|_{2,T,\bE}<\infty$, then the integral can be extended to $\bE$.
By the Burkholder-Davis-Gundy inequality for discontinuous martingales, for any $p \geq 1$
\begin{equation}
\label{BDG}E \sup_{t \leq T} |I^{\widehat{N}}(Y)(t,\bE)|^p \leq C_{p} E [I^{\widehat{N}}(Y)(\cdot,\bE)]_{T}^{p/2}.
\end{equation}

The previous inequality is not suitable for our purposes. A more convenient inequality can be obtained for {\em another} stochastic integral, constructed for $p \in [1,2]$ fixed, as suggested on page 293 of \cite{bie98}. More precisely, one can show that for any bounded simple process $Y$,
\begin{equation}
\label{moment-ineq2}E\sup_{t \leq T}|I^{\widehat{N}}(Y)(t,\bE)|^p \leq C_p E\int_0^T \int_{\bR^d} \int_{\bR \verb2\2 \{0\}}|Y(t,x,z)|^p \nu_{\alpha}(dz)dx dt=:[Y]_{p,T,\bE}^{p},
\end{equation}
where $C_p$ is the constant appearing in (\ref{BDG}) (see Lemma 8.22 of \cite{PZ07}).

By the usual procedure, the integral can be extended to the class of all $\cP \times \cB(\bR \verb2\2 \{0\})$-measurable processes $Y$ such that $[Y]_{p,T,\bE}<\infty$. The integral is defined as an element in the space $L^p(\Omega; D[0,T])$ and will be denoted by
$$I^{\widehat{N},p}(Y)(t,\bE)=\int_0^t \int_{\bR^d}\int_{\bR \verb2\2 \{0\}}Y(s,x,z)\widehat{N}(ds,dx,dz).$$
Its appealing feature is that it satisfies inequality (\ref{moment-ineq2}).

From now on, we fix $p \in [1,2]$. Based on (\ref{def-ZK-alpha-greater}), for any $B \in \cB_{b}(\bR^d)$, we let
$$I_{K}(X)(t,B)=\int_0^t \int_{B}X(s,x) Z_{K}(ds,dx)=\int_0^t \int_{B} \int_{\{ |z| \leq K\}} X(s,x)z \widehat{N}(ds,dx,dz),$$
for any predictable process $X=\{X(t,x);t \geq 0, x \in \bR^d\}$ for which the rightmost integral is well-defined. Letting $Y(t,x,z)=X(t,x)z 1_{\{0<|z| \leq K\}}$, we see that this is equivalent to saying that $p>\alpha$ and $X \in \cL_p$.
By (\ref{moment-ineq2}),
\begin{equation}
\label{p-moment-ineq2}E \sup_{t \leq T}|I_{K}(X)(t,B)|^p \leq C_{\alpha,p} K^{p-\alpha}E \int_0^T \int_{B}|X(s,x)|^p dx ds
\end{equation}
where $C_{\alpha,p}=C_p \alpha/(p-\alpha)$. If in addition, the process $X\in \cL_p$ satisfies
(\ref{p-finite-int-Rd})
then (\ref{p-moment-ineq2}) holds with $B$ replaced by $\cO$, for an arbitrary Borel set $\cO \subset \bR^d$.

Note that (\ref{p-moment-ineq2}) is the counterpart of (\ref{p-moment-ineq}) for the case $\alpha>1$. Together, these two inequalities will play a crucial role in Section \ref{nonlinear-eq-section}.

The following table summarizes all the conditions:
\begin{table}[ht]
\centering
\begin{tabular}{|c|c|c|} \hline
& $\alpha<1$ & $\alpha>1$ \\ \hline
$B$ is bounded & $X \in \cL_{\alpha}$ & $X \in \cL_{p}$ \\
 &  &  for some $p \in (\alpha,2]$ \\ \hline
$B=\cO$ is unbounded  & $X \in \cL_{\alpha}$ and  & $X \in \cL_{p}$ and \\
 & $X$ satisfies (\ref{finite-int-Rd})  & $X$ satisfies (\ref{p-finite-int-Rd}) \\
 &  & for some $p \in (\alpha,2]$\\
\hline
\end{tabular}
\caption{Conditions for $I_K(X)(t,B)$ to be well-defined}
\end{table}

\section{The main result}
\label{nonlinear-eq-section}

In this section, we state and prove the main result regarding the existence of a mild solution of equation (\ref{eq-noiseZ}). For this result, $\cO$ is a bounded domain in $\bR^d$. For any $t>0$,
we denote $$J_p(t)=\sup_{x \in \cO}\int_{\cO} G(t,x,y)^p dy.$$

\begin{theorem}
\label{main-th}
Let $\alpha \in (0,2)$, $\alpha \not=1$. Assume that for any $T>0$ 
\begin{equation}
\label{cond1}
\lim_{h \to 0}\int_0^T \int_{\cO}|G(t,x,y)-G(t+h,x,y)|^pdy dt = 0 \quad \forall x \in \cO,
\end{equation}
\begin{equation}
\label{cond2}
\lim_{|h| \to 0}\int_0^T \int_{\cO}|G(t,x,y)-G(t,x+h,y)|^pdy dt = 0 \quad \forall x \in \cO,
\end{equation}
\begin{equation}
\label{cond3}
\int_{0}^{T}J_p(t)dt<\infty,
\end{equation}
for some $p \in (\alpha,1)$ if $\alpha<1$, or for some $p \in (\alpha,2]$ if $\alpha>1$. 
Then equation (\ref{eq-noiseZ}) has a mild solution. Moreover, there exists a sequence $(\tau_K)_{K \geq 1}$ of stopping times with $\tau_K \uparrow \infty$ a.s. such that for any $T>0$ and $K \geq 1$,
$$\sup_{(t,x) \in [0,T] \times \cO}E(|u(t,x)|^p 1_{\{t \leq \tau_K \}})<\infty.$$
\end{theorem}

\begin{example} (Heat equation)
{\rm Let $L=\frac{\partial}{\partial t}-\frac{1}{2}\Delta$. Then $G(t,x,y) \leq \overline{G}(t,x-y)$ where $\overline{G}(t,x)$ is the fundamental solution of $Lu=0$ on $\bR^d$. Condition (\ref{cond3}) holds if $p<1+2/d$. If $\alpha<1$, this condition holds for any $p \in (\alpha,1)$. If $\alpha>1$, this condition holds for any $p \in (\alpha,1+2/d]$, as long as $\alpha$ satisfies (\ref{heat-cond}).
Conditions (\ref{cond1}) and (\ref{cond2}) hold by the continuity of the function $G$ in $t$ and $x$, by applying the dominated convergence theorem. To justify the application of this theorem, we use the trivial bound $(2\pi t)^{-dp/2}$ for both $G(t+h,x,y)^p$ and $G(t,x+h,y)^p$, which introduces the extra condition $dp<2$. Unfortunately, we could not find another argument for proving these two conditions. (In the case of the heat equation on $\bR^d$, Lemmas A.2 and A.3 of \cite{bie98} estimate the integrals appearing in (\ref{cond2}) and (\ref{cond1}), with $p=1$ in (\ref{cond1}). These arguments rely on the structure of $\overline{G}$ 
and cannot be used when $\cO$ is a bounded domain.)

}
\end{example}

\begin{example} (Parabolic equations)
{\rm Let $L=\frac{\partial }{\partial t}-\cL$ where $\cL$ is given by (\ref{def-cL-operator}). Assuming (\ref{parabolic-bound}), we see that
(\ref{cond3}) holds if $p<1+2/d$. The same comments as for the heat equation apply here as well. (Although in a different framework, a condition similar to (\ref{cond1}) was probably used in the proof of Theorem 12.11 of \cite{PZ07} (page 217) for the claim $\lim_{s \to t}E|J_3(X)(s)-J_3(X)(t)|_{L^p(\cO)}^p=0$. We could not see how to justify this claim, unless $dp<2$.)

}
\end{example}

\begin{example} (Heat equation with fractional power of the Laplacian)
{\rm Let $L=\frac{\partial}{\partial t}+(-\Delta)^{\gamma}$ for some $\gamma>0$.
By Lemma B.23 of \cite{PZ07}, if $\alpha>1$, then condition (\ref{cond3}) holds for any $p \in (\alpha,1+2\gamma/d)$, provided that $\alpha$ satisfies (\ref{heat-fractional-cond}). (This condition is the same as in Theorem 12.19 of \cite{PZ07}, which examines the same equation using the approach based on Hilbert-space valued solution.)
To verify condition (\ref{cond1}) and (\ref{cond2}), we use the continuity of $G$ in $t$ and $x$ and apply the dominated convergence theorem. To justify the application of this theorem, we use the trivial bound $C_{d,\gamma} t^{-dp/(2\gamma)}$ for both $G(t+h,x,y)^p$ and $G(t,x+h,y)^p$,  which introduces the extra condition $dp<2\gamma$.
This bound can be seen from (\ref{Green-function-fractional}), using the fact that $\cG(t,x,y) \leq \overline{\cG}(t,x-y)$ where $\cG$ and $\overline{\cG}$ are the fundamental solutions of $\frac{\partial u}{\partial t}-\Delta u=0$ on $\cO$, respectively $\bR^d$. (In the case of the same equation on $\bR^d$, elementary estimates for the time and space increments of $\overline{G}$ can be obtained directly from (\ref{Fourier-fractional}), as on p. 196 of \cite{azerad-mellouk07}.
These arguments cannot be used when $\cO$ is a bounded domain.)
}

\end{example}


The remaining part of this section is dedicated to the proof of Theorem \ref{main-th}.
The idea is to solve first the equation with the truncated noise $Z_{K}$ (yielding a mild solution $u_K$), and then identify a sequence $(\tau_K)_{K \geq 1}$ of stopping times with $\tau_{K} \uparrow \infty$ a.s. such that for any $t>0$, $x \in \cO$ and $L>K$, $u_{K}(t,x)=u_{L}(t,x)$ a.s. on the event $\{t \leq \tau_K\}$. The final step is to show that process $u$ defined by $u(t,x)=u_{K}(t,x)$ on $\{t \leq \tau_K\}$ is a mild solution of (\ref{eq-noiseZ}). A similar method can be found in Section 9.7 of \cite{PZ07} using an approach based on stochastic integration of operator-valued processes, with respect to Hilbert-space-valued processes, which is different from our approach.

Since $\sigma$ is a Lipschitz function, there exists a constant $C_{\sigma}>0$ such that:
\begin{equation}
\label{Lip1}
|\sigma(u)-\sigma(v)| \leq C_{\sigma}|u-v|, \quad \forall u,v \in \bR.
\end{equation}
In particular, letting $D_{\sigma}=C_{\sigma} \vee |\sigma(0)|$, we have:
\begin{equation}
\label{Lip2}
|\sigma(u)| \leq D_{\sigma}(1+|u|), \quad \forall u \in \bR.
\end{equation}

For the proof of Theorem \ref{main-th}, we need a specific construction of the Poisson random measure $N$, taken from \cite{PZ06}. We review briefly this construction.

Let $(\cO_k)_{k \geq 1}$ be a partition of $\bR^d$ with sets in $\cB_{b}(\bR^d)$ and $(U_j)_{j \geq 1}$ be a partition of $\bR \verb2\2 \{0\}$ such that $\nu_{\alpha}(U_j)<\infty$ for all $j \geq 1$. We may take $U_j=\Gamma_{j-1}$ for all $j \geq 1$. Let $(E_{i}^{j,k},X_{i}^{j,k},Z_{i}^{j,k})_{i,j,k \geq 1}$ be independent random variables defined on a probability space $(\Omega,\cF,P)$, such that
$$P(E_{i}^{j,k}>t)=e^{-\lambda_{j,k}t}, \quad P(X_{i}^{j,k} \in B)=\frac{|B \cap \cO_k|}{|\cO_k|}, \quad P(Z_{i}^{j,k} \in \Gamma)=\frac{|\Gamma \cap U_j|}{|U_j|},$$
where $\lambda_{j,k}=|\cO_k|\nu_{\alpha}(U_j)$.
Let $T_{i}^{j,k}=\sum_{l=1}^{i}E_{l}^{j,k}$ for all $i \geq 1$. Then
\begin{equation}
\label{construction-N}
N=\sum_{i,j,k \geq 1}\delta_{(T_i^{j,k},X_i^{j,k},Z_i^{j,k})}
\end{equation}
is a Poisson random measure on $\bR_{+} \times \bR^d \times (\bR \verb2\2 \{0\})$ with intensity $dtdx \nu_{\alpha}(dz)$.

This section is organized as follows. In Section \ref{eq-noise-ZK-section} we prove the existence of the solution of the equation with truncated noise $Z_K$. Sections \ref{alpha-less-section} and \ref{alpha-greater-section} contain the proof of Theorem \ref{main-th} when $\alpha<1$, respectively $\alpha>1$.

\subsection{The equation with truncated noise}
\label{eq-noise-ZK-section}

In this section, we fix $K>0$ and we consider the equation:
\begin{equation}
\label{eq-noise-ZK}Lu(t,x)=\sigma(u(t,x))\dot{Z}_{K}(t,x) \quad t>0,x \in \cO
\end{equation}
with zero initial conditions and Dirichlet boundary conditions. A mild solution of (\ref{eq-noise-ZK}) is a predictable process $u$ which satisfies (\ref{def-sol}) with $Z$ replaced by $Z_K$.
For the next result, $\cO$ can be a bounded domain in $\bR^d$ or $\cO=\bR^d$ (with no boundary conditions).

\begin{theorem}
\label{existence-th-ZK}
Under the assumptions of Theorem \ref{main-th}, equation (\ref{eq-noise-ZK}) has a unique mild solution $u=\{u(t,x);t \geq 0,x \in \cO\}$. For any $T>0$,
\begin{equation}
\label{moment-uK}
\sup_{(t,x) \in [0,T] \times \cO}E|u(t,x)|^p<\infty,
\end{equation}
and the map $(t,x) \mapsto u(t,x)$ is continuous from $[0,T] \times \cO$ into $L^p(\Omega)$.
\end{theorem}

\noindent {\bf Proof:} We use the same argument as in the proof of Theorem 13 of \cite{dalang99}, based on a Picard iteration scheme. We define $u_{0}(t,x)=0$ and
$$u_{n+1}(t,x)=\int_0^t \int_{\cO} G(t-s,x,y)\sigma(u_n(s,y))Z_{K}(ds,dy)$$
for any $n \geq 0$. We prove by induction on $n \geq 0$ that:
{\em (i)} $u_{n}(t,x)$ is well-defined;
{\em (ii)} $K_n(t):=\sup_{(t,x) \in [0,T] \times \cO}E|u_{n}(t,x)|^p<\infty$ for any $T>0$;\linebreak
{\em (iii)} $u_n(t,x)$ is $\cF_t$-measurable for any $t>0$ and $x \in \cO$;
{\em (iv)} the map $(t,x) \mapsto u_n(t,x)$ is continuous from $[0,T] \times \cO$ into $L^p(\Omega)$ for any $T>0$.

The statement is trivial for $n=0$. For the induction step, assume that the statement is true for $n$. 
By an extension to random fields of Theorem 30, Chapter IV of \cite{dellacherie-meyer75}, $u_n$ has a jointly measurable modification. Since this modification is $(\cF_t)_t$-adapted, (in the sense of {\em (iii)}), it has a predictable modification (using an extension of Proposition 3.21 of \cite{PZ07} to random fields). We work with this modification, that we call also $u_n$.

We prove that {\em (i)}-{\em (iv)} hold for $u_{n+1}$. To show {\em (i)}, it suffices to prove that $X_n \in \cL_{p}$, where $X_n(s,y)=1_{[0,t]}(s)G(t-s,x,y)\sigma(u_n(s,y))$. 
By (\ref{Lip2}) and (\ref{cond3}),
$$E\int_0^t \int_{\cO}|X_n(s,y)|^p dyds \leq D_{\sigma}^p 2^{p-1} (1+K_n(t)) \int_0^t J_p(t-s)ds<\infty.$$
In addition, if $\cO=\bR^d$, we have to prove that $X_n$ satisfies (\ref{finite-int-Rd}) if $\alpha<1$, or (\ref{p-finite-int-Rd}) if $\alpha>1$ (see Table 1). If $\alpha<1$, this follows as above, since $\alpha<p$ and hence $\sup_{(t,x) \in [0,T] \times \cO}E|u(t,x)|^{\alpha}<\infty$; the argument for $\alpha>1$ is similar.

\noindent Combined with the moment inequality (\ref{p-moment-ineq}) (or (\ref{p-moment-ineq2})), 
this proves {\em (ii)}, since:
\begin{equation}
\label{moment-bound-un}
E|u_{n+1}(t,x)|^p \leq C_{\alpha,p}K^{p-\alpha} D_{\sigma}^p 2^{p-1} (1+K_n(t)) \int_0^t J_p(t-s)ds,
\end{equation}
for any $x \in \cO$. Property {\em (iii)} follows by the construction of the integral $I_K$.

To prove {\em (iv)}, we first show the right continuity in $t$. Let $h>0$. Writing the interval $[0,t+h]$ as the union of $[0,t]$ and $(t,t+h]$, we obtain that $E|u_{n+1}(t+h,x)-u_{n+1}(t,x)|^p \leq 2^{p-1}(I_1(h)+I_2(h))$, where
\begin{eqnarray*}
I_1(h)&=& E\left|\int_0^t \int_{\cO}(G(t+h-s,x,y)-G(t-s,x,y)) \sigma(u_n(s,y)) Z_K(ds,dy) \right|^p \\
I_2(h)&=& E\left|\int_{t}^{t+h} \int_{\cO}G(t+h-s,x,y) \sigma(u_n(s,y)) Z_K(ds,dy) \right|^p.
\end{eqnarray*}
Using again (\ref{Lip2}) and the moment inequality (\ref{p-moment-ineq}) (or (\ref{p-moment-ineq2})), we obtain:
\begin{eqnarray*}
I_1(h) & \leq & D_{\sigma}^p 2^{p-1}(1+K_n(t)) \int_0^t \int_{\cO}|G(s+h,x,y)-G(s,x,y)|^p dyds \\
I_2(h) & \leq & D_{\sigma}^p 2^{p-1}(1+K_n(t)) \int_{0}^{h} \int_{\cO}G(s,x,y)^pdyds
\end{eqnarray*}
It follows that both $I_1(h)$ and $I_2(h)$ converge to $0$ as $h \to 0$, using (\ref{cond1}) for $I_1(h)$, respectively the Dominated Convergence Theorem and (\ref{cond3}) for $I_2(h)$. The left continuity in $t$ is similar, by writing the interval $[0,t-h]$ as the difference between $[0,t]$ and $(t-h,t]$ for $h>0$. For the continuity in $x$, similarly as above, we see that $E|u_{n+1}(t,x+h)-u_{n+1}(t,x)|^p$ is bounded by:
$$D_{\sigma}^p 2^{p-1} (1+K_n(t)) \int_0^t \int_{\cO}|G(s,x+h,y)-G(s,x,y)|^p dyds,$$
which converges to $0$ as $|h| \to 0$ due to (\ref{cond2}). This finishes the proof of {\em (iv)}.

We denote $M_n(t)=\sup_{x \in \cO}E|u_{n}(t,x)|^p$. Similarly to (\ref{moment-bound-un}), we have:
$$M_n(t) \leq C_1 \int_0^t (1+M_{n-1}(s))J_p(t-s)ds, \quad \forall n \geq 1,$$
where $C_1=C_{\alpha,p}K^{p-\alpha}D_{\sigma}^p 2^{p-1}$. By applying Lemma 15 of Erratum to \cite{dalang99} with $f_n=M_n$, $k_1=0$, $k_2=1$ and $g(s)=CJ_p(s)$, we obtain that:
\begin{equation}
\label{sup-M-finite}
\sup_{n \geq 0} \sup_{t \in [0,T]}M_n(t)<\infty \quad \mbox{for all} \ T>0.
\end{equation}

We now prove that $\{u_n(t,x)\}_n$ converges in $L^p(\Omega)$, uniformly in $(t,x) \in [0,T] \times \cO$. To see this, let $U_{n}(t)=\sup_{x \in \cO}E|u_{n+1}(t,x)-u_{n}(t,x)|^p$ for $n \geq 0$. Using the moment inequality (\ref{p-moment-ineq}) (or (\ref{p-moment-ineq2})) and (\ref{Lip1}), we have:
$$U_{n}(t) \leq C_2\int_0^t  U_{n-1}(s) J_p(t-s) ds$$
where $C_2=C_{\alpha,p}K^{p-\alpha} C_{\sigma}^p$. By Lemma 15 of Erratum to \cite{dalang99}, $\sum_{n \geq 0}U_n(t)^{1/p}$ converges uniformly on $[0,T]$. (Note that this lemma is valid for all $p>0$.)

We denote by $u(t,x)$ the limit of $u_{n}(t,x)$ in $L^p(\Omega)$. One can show that $u$ satisfies properties {\em (ii)}-{\em (iv)} listed above. So $u$ has a predictable modification. This modification is a solution of (\ref{eq-noise-ZK}). To prove uniqueness, let $v$ be another solution and denote $H(t)=\sup_{x \in \cO}E|u(t,x)-v(t,x)|^p$. Then
$$H(t) \leq C_2 \int_0^t H(s)J_p(t-s)ds.$$
Using (\ref{cond3}), it follows that $H(t)=0$ for all $t>0$. $\Box$

\subsection{Proof of Theorem \ref{main-th}: case $\alpha<1$}
\label{alpha-less-section}

In this case, for any $t>0$ and $B \in \cB_{b}(\bR^d)$, we have: (see (\ref{def-Z-alpha-less}))
$$Z(t,B)=\int_{[0,t] \times B \times (\bR \verb2\2 \{0\})}z N(ds,dx,dz).$$
The characteristic function of $Z(t,B)$ is given by:
$$E(e^{iuZ(t,B)})=\exp\left\{t|B| \int_{\bR \verb2\2 \{0\}} (e^{iuz}-1)\nu_{\alpha}(dz)\right\}, \quad \forall u \in \bR.$$
Note that $\{Z(t,B)\}_{t \geq 0}$ is {\em not} a compound Poisson process since  $\nu_{\alpha}$ is infinite.

We introduce the stopping times $(\tau_K)_{K \geq 1}$, as on page 239 of \cite{PZ06}:
$$\tau_{K}(B)=\inf\{t>0;|Z(t,B)-Z(t-,B)|>K\},$$
where $Z(t-,B)=\lim_{s \uparrow t}Z(s,B)$. Clearly, $\tau_{L}(B) \geq \tau_{K}(B)$ for all $L>K$.

We first investigate the relationship between $Z$ and $Z_K$ and the properties of $\tau_K(B)$.
Using construction (\ref{construction-N}) of $N$ and definition (\ref{def-ZK-alpha-less}) of $Z_K$, we have:
\begin{eqnarray*}
Z(t,B)&=&\sum_{i,j,k \geq 1}Z_{i}^{j,k} 1_{\{T_{i}^{j,k} \leq t\}} 1_{\{X_i^{j,k} \in B\}}=:\sum_{j,k \geq 1}Z^{j,k}(t,B)\\
Z_K(t,B)&=&\sum_{i,j,k \geq 1}Z_{i}^{j,k} 1_{\{|Z_{i}^{j,k}| \leq K\}} 1_{\{T_{i}^{j,k} \leq t\}} 1_{\{X_i^{j,k} \in B\}}.
\end{eqnarray*}

We observe that $\{Z^{j,k}(t,B)\}_{t \geq 0}$ is a compound Poisson process with
$$E(e^{iu Z^{j,k}(t,B)})=\exp\left\{t|\cO_k \cap B| \int_{U_j}(e^{iuz}-1)\nu_{\alpha}(dz) \right\} \quad \forall u \in \bR.$$

Note that $\tau_K(B)>T$ means that all the jumps of $\{Z(t,B)\}_{t \geq 0}$ in $[0,T]$ are smaller than $K$ in modulus, i.e.
$\{\tau_K(B)>T\}=\{\omega; |Z_i^{j,k}(\omega)| \leq K \ \mbox{for all} \linebreak i,j,k \geq 1 \ \mbox{for which} \ T_{i}^{j,k}(\omega) \leq T \ \mbox{and} \ X_{i}^{j,k}(\omega) \in B \}$. Hence, on $\{\tau_{K}(B)>T\}$,
$$Z([0,t] \times A)=Z_K([0,t] \times A)=Z_L([0,t] \times A),$$
for any $L>K$, $t \in [0,T]$, $A \in \cB_{b}(\bR^d)$ with $A \subset B$. Using an approximation argument and the construction of the integrals $I(X)$ and $I_K(X)$, it follows that for any $X \in \cL_{\alpha}$ and for any $L>K$, a.s. on $\{\tau_K(B)>T\}$, we have:
\begin{equation}
\label{Z-equal-ZK}
I(X)(T,B)=I_K(X)(T,B)=I_L(X)(T,B).
\end{equation}

The next result gives the probability of the event $\{\tau_K(B)>T\}$.

\begin{lemma}
\label{lim-tauK}
For any $T>0$ and $B \in \cB_{b}(\bR^d)$,
$$P(\tau_K(B)>T)=\exp(-T|B| K^{-\alpha}).$$
Consequently, $\lim_{K \to \infty}P(\tau_K(B)>T)=1$ and $\lim_{K \to \infty}\tau_K(B) = \infty$ a.s.
\end{lemma}

\noindent {\bf Proof:} Note that $\{\tau_{K}(B)>T\}=\bigcap_{j,k \geq 1}\{ \tau_{K}^{j,k}(B)>T\}$, where
$$\tau_{K}^{j,k}(B)=\inf\{t>0;|Z^{j,k}(t,B)-Z^{j,k}(t-,B)|>K\}$$
Since $\nu_{\alpha}(\{z;|z|>K\})=K^{-\alpha}$ and $(\tau_{K}^{j,k}(B))_{j,k \geq 1}$ are independent, it is enough to prove that for any $j,k \geq 1$,
\begin{equation}
\label{expo-tau-jk}P(\tau_K^{j,k}(B)>T)=\exp\left\{-T|B \cap \cO_k| \nu_{\alpha}(\{z;|z|>K\} \cap U_j) \right\}.
\end{equation}

Note that $\{\tau_{K}^{j,k}(B)>T\}=\{\omega;|Z_i^{j,k}(\omega)| \leq K \ \mbox{for all} \ i \ \mbox{for which} \ T_{i}^{j,k} \leq T \ \mbox{and}  \ X_{i}^{j,k} \in B\}$ and $(T_n^{j,k})_{n \geq 1}$ are the jump times of a Poisson process with intensity $\lambda_{j,k}$. Hence
\begin{eqnarray*}P(\tau_K^{j,k}(B)>T)
&=& \sum_{n \geq 0} \sum_{m=0}^{n} \sum_{I \subset \{1, \ldots,n\}, {\rm card}(I)=m}P(T_{n}^{j,k} \leq T<T_{n+1}^{j,k}) P(\bigcap_{i \in I}\{X_{i}^{j,k} \in B\})\\ & & P(\bigcap_{i \in I}\{|Z_i^{j,k}| \leq K\}) P(\bigcap_{i \in I^c}\{X_{i}^{j,k} \not \in B\} \\
&=& \sum_{n \geq 0}e^{-\lambda_{j,k}T} \frac{(\lambda_{j,k}T)^n}{n!} [1-P(X_{1}^{j,k} \in B)P(|Z_1^{j,k}|>K)]^n \\
&=& \exp\left\{-\lambda_{j,k}T P(X_{1}^{j,k} \in B)P(|Z_1^{j,k}|>K) \right\},
\end{eqnarray*}
which yields (\ref{expo-tau-jk}). 

To prove the last statement, let $A_k^{(n)}=\{\tau_{K}(B)>n\}$. Then $P(\overline{\lim}_{K}A_K^{(n)}) \geq \overline{\lim}_{K}P(A_{K}^{(n)})=1$ for any $n \geq 1$, and hence $P(\bigcap_{n \geq 1} \overline{\lim}_{K}A_K^{(n)})=1$. Hence, with probability $1$, for any $n$, there exists some $K_n$ such that $\tau_{K_n}>n$. Since $(\tau_K)_K$ is non-decreasing, this proves that $\tau_{K} \to \infty$ with probability $1$.
$\Box$

\begin{remark}
{\rm The construction of $\tau_K(B)$ given above is due to \cite{PZ06} (in the case of a symmetric measure $\nu_{\alpha}$). This construction relies on the fact that $B$ is a bounded set. Since $Z(t,\bR^d)$ (and consequently $\tau_{K}(\bR^d)$) is not well-defined, I could not see why this construction can also be used when $B=\bR^d$, as it is claimed in \cite{PZ06}. To avoid this difficulty, one could try to use an increasing sequence $(E_n)_n$ of sets in $\cB_{b}(\bR^d)$ with $\bigcup_{n}E_n=\bR^d$. Using (\ref{Z-equal-ZK}) with $B=E_n$ and letting $n \to \infty$, we obtain that $I(X)(t,\bR^d)=I_{K}(t,\bR^d)$ a.s. on $\{t \leq \tau_{K}\}$, where $\tau_K=\inf_{n \geq 1}\tau_{K}(E_n)$. But $P(\tau_K>t) \leq P(\underline{\lim}_{n} \{ \tau_{K}(E_n)>t\}) \leq \underline{\lim}_{n}P(\tau_{K}(E_n)>t)=\lim_{n}\exp(-t|E_n|K^{-\alpha})=0$ for any $t>0$, which means that $\tau_{K}=0$ a.s.
Finding a suitable sequence $(\tau_K)_K$ of stopping times which could be used in the case $\cO=\bR^d$ remains an open problem.
}
\end{remark}
\vspace{3mm}

In what follows, we denote $\tau_{K}=\tau_{K}(\cO)$. Let $u_{K}$ be the solution of equation (\ref{eq-noise-ZK}), whose existence is guaranteed by Theorem \ref{existence-th-ZK}.

\begin{lemma}
\label{uK-equal-uL}
Under the assumptions of Theorem \ref{main-th}, for any $t>0$, $x \in \cO$ and $L>K$,
$$u_{K}(t,x)=u_{L}(t,x) \quad \mbox{a.s. on} \ \{t \leq \tau_{K}\}.$$
\end{lemma}

\noindent {\bf Proof:} By the definition of $u_L$ and (\ref{Z-equal-ZK}),
\begin{eqnarray*}
u_{L}(t,x)&=&\int_0^t \int_{\cO}G(t-s,x,y)\sigma(u_{L}(s,y))Z_{L}(ds,dy) \\
&=&\int_0^t \int_{\cO}G(t-s,x,y)\sigma(u_{L}(s,y))Z_{K}(ds,dy)
\end{eqnarray*}
a.s. on the event $\{t \leq \tau_K\}$. Using the definition of $u_K$ and Proposition \ref{local-property-int} (Appendix \ref{appC-section}), we obtain that, with probability $1$,
\begin{eqnarray*}
(u_{K}(t,x)-u_{L}(t,x))1_{\{t \leq \tau_K\}}&=& 1_{\{t \leq \tau_K\}} \int_0^t \int_{\cO}G(t-s,x,y)(\sigma(u_{K}(s,y)) - \\
& & \quad \sigma(u_L(s,y)))1_{\{s \leq \tau_K\}} Z_K(ds,dy).
\end{eqnarray*}

Let $M(t)=\sup_{x \in \cO}E(|u_{K}(t,x)-u_{L}(t,x)|^p 1_{\{t \leq \tau_K\}})$.
Using the moment inequality (\ref{p-moment-ineq}) and the Lipschitz condition (\ref{Lip1}), we get:
$$M(t) \leq C \int_0^t J_p(t-s)M(s)ds,$$
where $C=C_{\alpha,p}K^{p-\alpha}C_{\sigma}^p$. Using (\ref{cond3}), it follows that $M(t)=0$ for all $t>0$. $\Box$

\vspace{3mm}

For any $t>0,x \in \cO$, let $\Omega_{t,x}=\bigcap_{L>K}\{t \leq \tau_{K}(t), u_K(t,x) \not=u_{L}(t,x)\}$,
where $L$ and $K$ are positive integers.
Let $\Omega_{t,x}^*=\Omega_{t,x} \cap \{\lim_{K \to \infty}\tau_K=\infty\}$.
 By Lemma \ref{lim-tauK} and Lemma \ref{uK-equal-uL}, $P(\Omega_{t,x}^*)=1$.


\vspace{3mm}

The next result concludes the proof of Theorem \ref{main-th}.

\begin{proposition}
\label{u-is-solution}
Under the assumptions of Theorem \ref{main-th},
the process $u=\{u(t,x);t \geq  0,x \in \cO\}$ defined by:
\begin{eqnarray*}
u(\omega,t,x)&=&u_{K}(\omega,t,x) \quad if \ \omega \in \Omega_{t,x}^* \ \mbox{and} \ t \leq \tau_K(\omega)\\
u(\omega,t,x)&=& 0 \quad \quad \quad \quad \quad if \ \omega \not\in \Omega_{t,x}^*
\end{eqnarray*}
is a mild solution of equation (\ref{eq-noiseZ}).
\end{proposition}

\noindent {\bf Proof:} We first prove that $u$ is predictable. Note that
$$u(t,x)=\lim_{K \to \infty}(u_{K}(t,x) 1_{\{t \leq \tau_K\}})1_{\Omega_{t,x}^*}.$$
The process $X(\omega,t,x)=1_{\{t \leq \tau_K\}}(\omega)$ is clearly predictable, being in the class $\cC$ defined in Remark \ref{leftcont-predictable}.
By the definition of $\Omega_{t,x}$, since $u_K,u_L$ are predictable, it follows that $(\omega,t,x) \mapsto 1_{\Omega_{t,x}^*}(\omega)$ is $\cP$-measurable. Hence, $u$ is predictable.

We now prove that $u$ satisfies (\ref{def-sol}). Let $t>0$ and $x \in \cO$ be arbitrary. Using (\ref{Z-equal-ZK}) and Proposition \ref{local-property-int} (Appendix \ref{appC-section}), with probability $1$, we have:
\begin{eqnarray*}
1_{\{t \leq \tau_K\}} u(t,x)  &=& 1_{\{t \leq \tau_K\}} u_{K}(t,x) \\
 &=& 1_{\{t \leq \tau_K\}} \int_0^t \int_{\cO} G(t-s,x,y) \sigma(u_{K}(s,y)) Z_K(ds,dy) \\
 &=& 1_{\{t \leq \tau_K\}} \int_0^t \int_{\cO} G(t-s,x,y) \sigma(u_{K}(s,y)) Z(ds,dy) \\
 &=&  1_{\{t \leq \tau_K\}} \int_0^t \int_{\cO} G(t-s,x,y) \sigma(u_{K}(s,y)) 1_{\{s \leq \tau_K\}}Z(ds,dy)\\
 &=& 1_{\{t \leq \tau_K\}} \int_0^t \int_{\cO} G(t-s,x,y) \sigma(u(s,y)) 1_{\{s \leq \tau_K\}}Z(ds,dy) \\
 &=& 1_{\{t \leq \tau_K\}} \int_0^t \int_{\cO} G(t-s,x,y) \sigma(u(s,y)) Z(ds,dy).
\end{eqnarray*}
For the second last equality, we used the fact that processes $X(s,y)=1_{[0,t]}(s)G(t-s,x,y) \sigma(u_{K}(s,y)) 1_{\{s \leq \tau_K\}}$ and $Y(s,y)=1_{[0,t]}(s)G(t-s,x,y) \sigma(u(s,y)) \linebreak 1_{\{s \leq \tau_K\}}$ are modifications of each other (i.e. $X(s,y)=Y(s,y)$ a.s. for all $s>0,y \in \cO$), and hence, $[X-Y]_{\alpha,t,\cO}=0$ and $I(X)(t,\cO)=I(Y)(t,\cO)$ a.s. The conclusion follows letting $K \to \infty$, since $\tau_K \to \infty$ a.s. $\Box$

\subsection{Proof of Theorem \ref{main-th}: case $\alpha>1$}
\label{alpha-greater-section}

In this case, for any $t>0$ and $B \in \cB_b(\bR^d)$, we have: (see (\ref{def-Z-alpha-greater}))
$$Z(t,B)=
\int_{[0,t] \times B \times (\bR \verb2\2 \{0\})}z \widehat{N}(ds,dx,dz).$$

To introduce the stopping times $(\tau_K)_{K \geq 1}$ we use the same idea as in Section 9.7 of \cite{PZ07}.

Let $M(t,B)=\sum_{j \geq 1}(L_j(t,B)-EL_j(t,B))$ and $P(t,B)=L_0(t,B)$, where $L_j(t,B)=L_j([0,t] \times B)$ was defined in Section \ref{def-noise-section}. Note that
$\{M(t,B)\}_{t \geq 0}$ is a zero-mean square-integrable martingale
and $\{P(t,B)\}_{t \geq 0}$ is a compound Poisson process with
$E[P(t,B)]=t|B|\mu$ where $\mu=\int_{|z|>1}z \nu_{\alpha}(dz)=\beta \frac{\alpha}{\alpha-1}$.
With this notation,
\begin{equation}
\label{Z-split}Z(t,B)=M(t,B)+P(t,B)-t|B|\mu.
\end{equation}

We let $M_{K}(t,B)=P_K(t,B)-E[P_K(t,B)]=P_{K}(t,B)-t|B|\mu_{K}$, where
$$P_{K}(t,B)=\int_{[0,t] \times B \times (\bR \verb2\2 \{0\})}z 1_{\{1 <|z| \leq K \}}N(ds,dx,dz)$$ and $\mu_{K}=\int_{1<|z| \leq K}z \nu_{\alpha}(dz)$.
Recalling definition (\ref{def-ZK-alpha-greater}) of $Z_K$, it follows that:
\begin{equation}
\label{ZK-split}
Z_{K}(t,B)=
M(t,B)+P_{K}(t,B)-t|B|\mu_{K}.
\end{equation}

For any $K>0$, we let
$$\tau_K(B)=\inf\{t>0;|P(t,B)-P(t-,B)|>K \},$$
where $P(t-,B)=\lim_{s \uparrow t}P(s,B)$.

Lemma \ref{lim-tauK} holds again, but its proof is simpler than in the case $\alpha<1$, since $\{P(t,B)\}_{t \geq 0}$ is a compound Poisson process. By (\ref{construction-N}), 
\begin{eqnarray*}
P(t,B)&=&\sum_{i,j,k \geq 1}Z_{i}^{j,k} 1_{\{|Z_i^{j,k}|>1\}} 1_{\{T_i^{j,k} \leq t \}} 1_{\{X_{i}^{j,k} \in B \}} \\
P_K(t,B)&=& \sum_{i,j,k \geq 1}Z_{i}^{j,k} 1_{\{1<|Z_i^{j,k}| \leq K\}} 1_{\{T_i^{j,k} \leq t \}} 1_{\{X_{i}^{j,k} \in B \}}.
\end{eqnarray*}
Hence, on $\{\tau_K(B)>T\}$, for any $L>K$, $t \in [0,T]$, $A \in \cB_{b}(\bR^d)$ with $A \subset B$,
$$P([0,t] \times A)=P_{K}([0,t] \times A)=P_{L}([0,t] \times A).$$

Let $b_{K}=\mu-\mu_{K}=\int_{|z|>K}z \nu_{\alpha}(dz)$. Using (\ref{Z-split}) and (\ref{ZK-split}), it follows that:
$$Z([0,t] \times A)=Z_K([0,t] \times A)-t|A| b_K=Z_L([0,t] \times A)-t|A|b_L$$
for any $L>K$, $t \in [0,T]$, $A \in \cB_{b}(\bR^d)$ with $A \subset B$. Let $p \in (\alpha,2]$ be fixed. Using an approximation argument and the construction of the integrals $I(X)$ and $I_K(X)$, it follows that for any $X \in \cL_{\alpha}$ and for any $L>K$, a.s. on $\{\tau_K(B)>T\}$, we have:
\begin{eqnarray}
\label{Z-equal-ZK-part2}
I(X)(T,B)&=&I_K(X)(T,B)-b_K \int_0^T \int_{\cO} X(s,y)dyds \\
\nonumber
&=& I_L(X)(T,B) -b_L \int_0^T \int_{\cO} X(s,y)dyds.
\end{eqnarray}

We denote $\tau_K=\tau_{K}(\cO)$. We consider the following equation:
\begin{equation}
\label{eq-noise-ZK-part2}Lu(t,x)=\sigma(u(t,x))\dot{Z}_K(t,x)-b_{K} \sigma(u(t,x)), \quad t>0,x \in \cO
\end{equation}
with zero initial conditions and Dirichlet boundary conditions. A mild solution of (\ref{eq-noise-ZK-part2}) is a predictable process $u$ which satisfies:
\begin{eqnarray*}
u(t,x)&=&\int_0^t \int_{\cO}G(t-s,x,y)\sigma(u(s,y))Z_{K}(ds,dy)-\\
& & b_K \int_0^t \int_{\cO}G(t-s,x,y)\sigma(u(s,y)) dyds \quad {\rm a.s.}
\end{eqnarray*}
for any $t>0,x \in \cO$. The existence and uniqueness of a mild solution of (\ref{eq-noise-ZK-part2}) can be proved similarly to Theorem \ref{existence-th-ZK}. We omit these details. We denote this solution by $v_K$.

\begin{lemma}
\label{vK-equal-vL}
Under the assumptions of Theorem \ref{main-th}, for any $t>0,x \in \cO$ and $L>K$,
$$v_K(t,x)=v_{L}(t,x) \quad a.s. \ on \ \{t \leq \tau_K\}.$$
\end{lemma}

\noindent {\bf Proof:} By the definition of $v_L$ and (\ref{Z-equal-ZK-part2}), a.s. on the event $\{t \leq \tau_{K}\}$, $v_{L}(t,x)$ is equal to
$$\int_0^t \int_{\cO}G(t-s,x,y)\sigma(v_{L}(s,y))Z_{L}(ds,dy)-b_{L} \int_0^t \int_{\cO}G(t-s,x,y)\sigma(v_{L}(s,y))dyds=$$
$$\int_0^t \int_{\cO}G(t-s,x,y)\sigma(v_{L}(s,y))Z_{K}(ds,dy)-b_{K} \int_0^t \int_{\cO}G(t-s,x,y)\sigma(v_{L}(s,y))dyds.$$

Using the definition of $v_K$ and Proposition \ref{local-property-int} (Appendix \ref{appC-section}), we obtain that, with probability $1$,
\begin{eqnarray*}
\lefteqn{(v_{K}(t,x)-v_{L}(t,x))1_{\{t \leq \tau_K\}} =1_{\{t \leq \tau_K\}} \left( \int_0^t \int_{\cO} G(t-s,x,y) (\sigma(v_{K}(s,y))-\sigma(v_{L}(s,y)) \right. } \\
& & \left. 1_{\{s \leq \tau_K\}} Z_{K}(ds,dy)-\int_0^t \int_{\cO} G(t-s,x,y)(\sigma(v_{K}(s,y))-\sigma(v_{L}(s,y)) 1_{\{s \leq \tau_K\}}dyds \right).
\end{eqnarray*}

Letting $M(t)=\sup_{x \in \cO}E(|v_K(t,x)-v_L(t,x)|^p 1_{\{t \leq \tau_K\}})$, we see that $M(t) \leq 2^{p-1}(E|A(t,x)|^p+E|B(t,x)|^p)$ where
\begin{eqnarray*}
A(t,x)&=& \int_0^t \int_{\cO} G(t-s,x,y)(\sigma(v_{K}(s,y))-\sigma(v_{L}(s,y)) 1_{\{s \leq \tau_K\}} Z_{K}(ds,dy)  \\
B(t,x)&=&  \int_0^t \int_{\cO} G(t-s,x,y)(\sigma(v_{K}(s,y))-\sigma(v_{L}(s,y)) 1_{\{s \leq \tau_K\}} dyds.
\end{eqnarray*}

We estimate separately the two terms. For the first term, we use the moment inequality (\ref{p-moment-ineq2}) and the Lipschitz condition (\ref{Lip1}). We get:
\begin{equation}
\label{sup-mean-A}
\sup_{x \in \cO}E|A(t,x)|^p \leq C \int_0^t J_p(t-s)M(s)ds,
\end{equation}
where $C=C_{\alpha,p}K^{p-\alpha}C_{\sigma}^p$. For the second term, we use H\"older's inequality $|\int fg d\mu| \leq (\int |f|^p d\mu )^{1/p} (\int |g|^q d\mu )^{1/q}$ with
$f(s,y)=G(t-s,x,y)^{1/p} (\sigma(v_{K}(s,y))-\sigma(v_{L}(s,y)) 1_{\{s \leq \tau_K\}}$ and $g(s,y)=G(t-s,x,y)^{1/q}$, where $p^{-1}+q^{-1}=1$. Hence,
$$|B(t,x)|^p \leq C_{\sigma}^p K_t^{p/q}\int_0^t \int_{\cO}G(t-s,x,y)|v_K(s,y)-v_L(s,y)|^p 1_{\{s \leq \tau_K\}}dyds,$$
where $K_t=\int_0^t J_1(s)ds<\infty$. (Since $\cO$ is a bounded set, $J_1(s) \leq C J_p(s)^{1/p}$ where $C$ is a constant depending on $|\cO|$ and $p$. Since $p>1$, $\int_0^t J_p(s)^{1/p}ds \leq c_t (\int_0^t J_p(s)ds)^{1/p}<\infty$ by (\ref{cond3}). This shows that $K_t<\infty$.) Therefore,
\begin{equation}
\label{sup-mean-B}
\sup_{x \in \cO}E|B(t,x)|^p \leq  C_t \int_0^t J_1(t-s) M(s)ds,
\end{equation}
where $C_t=C_{\sigma}^p K_t^{p/q}$.
From (\ref{sup-mean-A}) and (\ref{sup-mean-B}), we obtain that:
$$M(t) \leq C_t' \int_0^t (J_p(t-s)+J_1(t-s))M(s)ds$$
where $C_t'=2^{p-1} (C \vee C_t)$. This implies that $M(t)=0$ for all $t>0$. $\Box$

\vspace{3mm}

For any $t>0$ and $x \in \cO$, we let $\Omega_{t,x}=\bigcap_{L>K}\{t \leq \tau_K,v_K(t,x) \not=v_L(t,x)\}$ where $K$ and $L$ are positive integers, and $\Omega_{t,x}^*=\Omega_{t,x} \cap\{\lim_{K \to \infty}\tau_K=\infty \}$. By Lemma \ref{vK-equal-vL}, $P(\Omega_{t,x}^*)=1$.

\begin{proposition}
Under the assumptions of Theorem \ref{main-th}, the process $u=\{u(t,x); t \geq 0, x \in \cO\}$ defined by:
\begin{eqnarray*}
u(\omega,t,x)&=&v_{K}(\omega,t,x) \quad if \ \omega \in \Omega_{t,x}^* \ \mbox{and} \ t \leq \tau_K(\omega)\\
u(\omega,t,x)&=& 0 \quad \quad \quad \quad \quad if \ \omega \not\in \Omega_{t,x}^*
\end{eqnarray*}
is a mild solution of equation (\ref{eq-noiseZ}).
\end{proposition}

\noindent {\bf Proof:} We proceed as in the proof of Proposition \ref{u-is-solution}. In this case, with probability $1$, we have:
\begin{eqnarray*}
1_{\{t \leq \tau_K\}}u(t,x)&=&1_{\{t \leq \tau_K\}} \left(\int_0^{t} \int_{\cO}G(t-s,x,y)\sigma(u(s,y))Z(ds,dy)- \right. \\
& & \left. b_{K} \int_0^t \int_{\cO} G(t-s,x,y)\sigma(u(s,y))dyds \right).
\end{eqnarray*}
The conclusion follows letting $K \to \infty$, since $\tau_K \to \infty$ a.s. and $b_{K} \to 0$. $\Box$

\appendix
\section{Some auxiliary results}
\label{appA-section}

The following result is used in the proof of Theorem \ref{max-ineq}.
\begin{lemma}
\label{lemmaA1}
If $X$ has a $S_{\alpha}(\sigma,\beta,0)$ distribution then
$$\lambda^{\alpha}P(|X|>\lambda) \leq c_{\alpha}^*\sigma^{\alpha} \quad \mbox{for all} \ \lambda>0,$$
where $c_{\alpha}^*>0$ is a constant depending only on $\alpha$.
\end{lemma}

\noindent {\bf Proof:} {\em Step 1.} We first prove the result for $\sigma=1$. We treat only the right tail, the left tail being similar. We denote $X$ by $X_{\beta}$ to emphasize the dependence on $\beta$. By Property 1.2.15 of \cite{ST94}, $\lim_{\lambda \to \infty}\lambda^{\alpha}P(X_{\beta}>\lambda)=C_{\alpha}\frac{1+\beta}{2}$, where $C_{\alpha}=(\int_0^{\infty}x^{-\alpha}\sin x dx)^{-1}$. We use the fact that for any $\beta \in [-1,1]$,
$$P(X_{\beta} > \lambda) \leq P(X_{1} > \lambda) \quad \mbox{for all} \ \lambda>\lambda_{\alpha}$$ 
for some $\lambda_{\alpha}>0$ (see Property 1.2.14 of \cite{ST94} or Section 1.5 of \cite{nolan13}). 
Since $\lim_{\lambda \to \infty}\lambda^{\alpha}P(X_1>\lambda)=C_{\alpha}$, there exists $\lambda_{\alpha}^*>\lambda_{\alpha}$ such that $$\lambda^{\alpha}P(X_1>\lambda) <2C_{\alpha}\quad \mbox{for all} \ \lambda>\lambda_{\alpha}^*.$$ 
It follows that $\lambda^{\alpha}P(X_{\beta}>\lambda)<2C_{\alpha}$ for all $\lambda >\lambda_{\alpha}^*$ and $\beta \in [-1,1]$. Clearly, for all $\lambda \in (0,\lambda_{\alpha}^*]$ and $\beta \in [-1,1]$, $\lambda^{\alpha}P(X_{\beta}>\lambda) \leq \lambda^{\alpha} \leq (\lambda_\alpha^*)^{\alpha}$.

{\em Step 2.} We now consider the general case. Since $X/\sigma$ has a $S_{\alpha}(1,\beta,0)$ distribution, by Step 1, it follows that $\lambda^{\alpha}P(|X|>\sigma \lambda) \leq c_{\alpha}^*$ for any $\lambda>0$. The conclusion follows multiplying by $\sigma^{\alpha}$. $\Box$

\vspace{3mm}

In the proof of Theorem \ref{max-ineq} and Lemma \ref{lemmaA2} below, we use the following remark, due to Adam Jakubowski (personal communication).

\begin{remark}
\label{Adam-remark}
{\rm Let $X$ be a random variable such that $P(|X|>\lambda) \leq K \lambda^{-\alpha}$ for all $\lambda>0$, for some $K>0$ and $\alpha \in (0,2)$. Then, for any $A>0$,
$$E(|X| 1_{\{|X| \leq A\}}) \leq \int_0^A P(|X|>t)dt \leq K \frac{1}{1-\alpha}A^{1-\alpha} \quad \mbox{if} \ \alpha <1,$$
$$E(|X| 1_{\{|X| > A\}}) \leq \int_{A}^{\infty} P(|X|>t)dt +AP(|X|>A) \leq K\frac{\alpha}{\alpha-1}A^{1-\alpha} \quad \mbox{if} \ \alpha >1,$$
$$E(X^2 1_{\{|X| \leq A\}}) \leq 2 \int_{0}^{A} t P(|X|>t)dt \leq K\frac{2}{2-\alpha}A^{2-\alpha} \quad \mbox{for any} \ \alpha \in (0,2).$$
}
\end{remark}

The next result is a generalization of Lemma 2.1 of \cite{GM82} to the case of non-symmetric random variables. This result is used in the proof of Lemma \ref{tail-ZKB-lemma} and Proposition \ref{max-ineq-alpha-less-1}.

\begin{lemma}
\label{lemmaA2}
Let $(\eta_k)_{k \geq 1}$ be independent random variables such that
\begin{equation}
\label{tail-prob-term}
\sup_{\lambda>0}\lambda^{\alpha}P(|\eta_k|>\lambda) \leq K \quad \forall k \geq 1
\end{equation}
for some $K>0$ and $\alpha \in (0,2)$. If $\alpha>1$, we assume that $E(\eta_k)=0$ for all $k$, and if $\alpha=1$, we assume that $\eta_k$ has a symmetric distribution for all $k$. Then for any sequence $(a_k)_{k \geq 1}$ of real numbers, we have:
\begin{equation}
\label{tail-prob-sum}
\sup_{\lambda>0}\lambda^{\alpha}P(|\sum_{k \geq 1}a_k \eta_k|>\lambda) \leq r_{\alpha}K \sum_{k \geq 1}|a_k|^{\alpha}
\end{equation}
where $r_{\alpha}>0$ is a constant depending only on $\alpha$.
\end{lemma}

\noindent {\bf Proof:} We consider the intersection of the event on the left-hand side of (\ref{tail-prob-sum}) with the event $\{\sup_{k \geq 1}|a_k \eta_k|>\lambda\}$ and its complement. Hence,
$$P(|\sum_{k \geq 1}a_k \eta_k|>\lambda) \leq \sum_{k \geq 1}P(|a_k\eta_k|>\lambda)+P(|\sum_{k \geq 1}a_k \eta_k 1_{\{|a_k\eta_k| \leq \lambda\}}|>\lambda)=:I+II.$$
Using (\ref{tail-prob-term}), we have
$I \leq K \lambda^{-\alpha}\sum_{k \geq 1}|a_k|^{\alpha}$. To treat $II$, we consider 3 cases.

{\em Case 1.} $\alpha<1$. By Markov's inequality and Remark \ref{Adam-remark}, we have:
$$II \leq \frac{1}{\lambda}\sum_{k \geq 1}|a_k| E(|\eta_k| 1_{\{|a_k \eta_k| \leq \lambda\}}) \leq K \frac{1}{1-\alpha} \lambda^{-\alpha}\sum_{k \geq 1}|a_k|^{\alpha}.$$

{\em Case 2.} $\alpha>1$. Let $X=\sum_{k \geq 1}a_k \eta_k 1_{\{|a_k \eta_k| \leq \lambda\}}$. Since $E(\sum_{k \geq 1}a_k \eta_k)=0$,
$$|E(X)|=|E(\sum_{k \geq 1}a_k \eta_k 1_{\{|a_k \eta_k| > \lambda\}})| \leq \sum_{k \geq 1}|a_k| E(|\eta_k| 1_{\{|a_k \eta_k|>\lambda\}}) \leq \frac{K\alpha}{\alpha-1}\lambda^{1-\alpha}\sum_{k \geq 1}|a_k|^{\alpha},$$
where we used Remark \ref{Adam-remark} for the last inequality.
From here, we infer that
$$|E(X)|<\frac{\lambda}{2} \quad \mbox{for any} \ \lambda>\lambda_{\alpha},$$ where $\lambda_{\alpha}^{\alpha}=2K \frac{\alpha}{\alpha-1}\sum_{k \geq 1}|a_k|^{\alpha}$. By Chebyshev's inequality, for any $\lambda>\lambda_{\alpha}$,
\begin{eqnarray*}
II &=& P(|X|>\lambda) \leq  P(|X-E(X)|>\lambda-|E(X)|) \leq \frac{4}{\lambda^2} E|X-E(X)|^2 \\
& \leq & \frac{4}{\lambda^2}\sum_{k \geq 1}a_k^2 E(\eta_k^2 1_{\{|a_k \eta_k| \leq \lambda \}}) \leq \frac{8K}{2-\alpha}\lambda^{-\alpha}\sum_{k \geq 1}|a_k|^{\alpha},
\end{eqnarray*}
using Remark \ref{Adam-remark} for the last inequality. On the other hand, if $\lambda \in (0,\lambda_{\alpha}]$,
$$II=P(|X|>\lambda) \leq 1 \leq \lambda_{\alpha}^{\alpha}\lambda^{-\alpha}=2K \frac{\alpha}{\alpha-1}\lambda^{-\alpha}\sum_{k \geq 1}|a_k|^{\alpha}.$$

{\em Case 3.} $\alpha=1$. Since $\eta_k$ has a symmetric distribution, we can use the original argument of \cite{GM82}. $\Box$

\section{Fractional power of the Laplacian}
\label{appB-section}

Let $\overline{G}(t,x)$ be the fundamental solution of $\frac{\partial u}{\partial t}+(-\Delta)^{\gamma}u=0$ on $\bR^d$, $\gamma>0$.

\begin{lemma}
\label{fract-lemma}
For any $p>1$, there exist some constants $c_1,c_2>0$  depending on $d,p,\gamma$ such that
$$c_1 t^{-\frac{d}{2\gamma}(p-1)} \leq \int_{\bR^d} \overline{G}(t,x)^{p}dx \leq c_2 t^{-\frac{d}{2\gamma}(p-1)}.$$
\end{lemma}

\noindent {\bf Proof:} The upper bound is given by Lemma B.23 of \cite{PZ07}. For the lower bound, we use the scaling property of the functions $(g_{t,\gamma})_{t>0}$. We have:
\begin{eqnarray*}
\overline{G}(t,x)&=&\int_0^{\infty} \frac{1}{(4\pi t^{1/\gamma}r)^{d/2}} \exp\left(-\frac{|x|^2}{4t^{1/\gamma}r} \right)g_{1,\gamma}(r)dr\\
& \geq & \int_1^{\infty} \frac{1}{(4\pi t^{1/\gamma}r)^{d/2}} \exp\left(-\frac{|x|^2}{4t^{1/\gamma}r} \right)g_{1,\gamma}(r)dr\\
& \geq & \frac{1}{(4\pi t^{1/\gamma})^{d/2}} \exp\left(-\frac{|x|^2}{4t^{1/\gamma}} \right)C_{d,\gamma} \quad \mbox{with} \quad C_{d,\gamma}:=\int_1^{\infty}r^{-d/2} g_{1,\gamma}(r)dr<\infty,
\end{eqnarray*}
and hence
$$\int_{\bR^d}\overline{G}(t,x)^p dx \geq c_{d,\gamma,p}' t^{-\frac{dp}{2\gamma}}\int_{\bR^d}\exp\left(-\frac{p|x|^2}{4t^{1/\gamma}}\right)dx=c_{d,p,\gamma}
t^{-\frac{d}{2\gamma}(p-1)}.$$
$\Box$

\section{A local property of the integral}
\label{appC-section}

The following result is the analogue of Proposition 8.11 of \cite{PZ07}.

\begin{proposition}
\label{local-property-int}
Let $T>0$ and $\cO \subset \bR^d$ be a Borel set. Let $X=\{X(t,x);\linebreak t \geq 0,x \in \bR^d\}$ be a predictable process such that $X \in \cL_{\alpha}$ if $\alpha<1$, or $X \in \cL_{p}$ for some $p \in (\alpha,2]$ if $\alpha>1$. If $\cO$ is unbounded, assume in addition that $X$ satisfies (\ref{finite-int-Rd}) if $\alpha<1$, or $X$ satisfies (\ref{p-finite-int-Rd}) for some $p \in (\alpha,2)$, if $\alpha>1$. Suppose that there exists an event $A \in \cF_T$ such that
\begin{equation}
\label{X-is-zero}
X(\omega,t,x)=0 \quad for all \ \omega \in A, t \in [0,T], x \in \cO.
\end{equation}
Then for any $K>0$, $I(X)(T,\cO)=I_{K}(X)(T,\cO)=0$ a.s. on $A$.
\end{proposition}

\noindent {\bf Proof:} We only prove the result for $I(X)$, the proof for $I_{K}(X)$ being the same. Moreover, we include only the argument for $\alpha<1$; the case $\alpha>1$ is similar. The idea is to reduce the argument to the case when $X$ is a simple process, as in the proof Proposition of 8.11 of \cite{PZ07}.

{\em Step 1.} We show that the proof can be reduced to the case of a bounded set $\cO$. Let $X_n(t,x)=X(t,x)1_{\cO_n}(x)$ where $\cO_n=\cO \cap E_n$ and $(E_n)_{n}$ ia an increasing sequence of sets in $\cB_{b}(\bR^d)$ such that $\bigcup_n E_n=\bR^d$. Then $X_n \in \cL_{\alpha}$ satisfies (\ref{X-is-zero}). By the dominated convergence theorem, $$E\int_0^T \int_{\cO}|X_n(t,x)-X(t,x)|^{\alpha} \to 0.$$ By the construction of the integral, $I(X_{n_k})(T,\cO) \to I(X)(T,\cO)$ a.s. for a subsequence $\{n_k\}$. It suffices to show that $I(X_{n})(T,\cO)=0$ a.s. on $A$ for all $n$. But $I(X_{n})(T,\cO)=I(X_{n})(T,\cO_n)$ and $\cO_n$ is bounded.

{\em Step 2.} We show that the proof can be reduced to the case of a bounded processes. For this, let $X_n(t,x)=X(t,x)1_{\{|X(t,x)|\leq n\}}$. Clearly, $X_n \in \cL_{\alpha}$ is bounded and satisfies (\ref{X-is-zero}) for all $n$. By the dominated convergence theorem, $[X_n-X]_{\alpha} \to 0$,
and hence $I(X_{n_k})(T,\cO) \to I(X)(T,\cO)$ a.s. for a subsequence $\{n_k\}$. It suffices to show that $I(X_{n})(T,\cO)=0$ a.s. on $A$ for all $n$.

{\em Step 3.} We show that the proof can be reduced to the case of bounded continuous processes. Assume that $X \in \cL_{\alpha}$ is bounded and satisfies (\ref{X-is-zero}). For any $t>0$ and $x \in \bR^d$, we define
$$X_{n}(t,x)=n^{d+1}\int_{(t-1/n) \vee 0}^{t} \int_{(x-1/n,x] \cap \cO}X(s,y)dyds,$$
where $(a,b]=\{y \in \bR^d; a_i<y_i \leq b_i \ \mbox{for all} \ i=1,\ldots,d\}$.
Clearly, $X_n$ is bounded and satisfies (\ref{X-is-zero}).
We prove that $X_n \in \cL_{\alpha}$. Since $X_n$ is bounded, $[X_n]_{\alpha}<\infty$. To prove that $X_n$ is predictable, we consider
$$F(t,x)=\int_0^t \int_{(0,x] \cap \cO}X(s,y)dyds.$$
Since $X$ is predictable, it is progressively measurable, i.e. for any $t>0$, the map $(\omega,s,x) \mapsto X(\omega,s,x)$ from $\Omega \times [0,t] \times \bR^d$ to $\bR$ is $\cF_t \times \cB([0,t]) \times \cB(\bR^d)$-measurable. Hence,
$F(t,\cdot)$ is $\cF_t \times \cB(\bR^d)$-measurable for any $t>0$. Since the map $t \mapsto F(\omega,t,x)$ is left-continuous for any $\omega \in \Omega,x \in \bR^d$, it follows that $F$ is predictable, being in the class $\cC$ defined in Remark \ref{leftcont-predictable}. Hence, $X_n$ is predictable, being a sum of $2^{d+1}$ terms involving $F$.

Since $F$ is continuous in $(t,x)$, $X_n$ is continuous in $(t,x)$.
By Lebesque differentiation theorem in $\bR^{d+1}$, $X_n(\omega,t,x) \to X(\omega,t,x)$ for any $\omega \in \Omega, t>0,x \in \cO$. By the bounded convergence theorem, $[X_n-X]_{\alpha} \to 0$. Hence $I(X_{n_k})(T,\cO) \to I(X)(T,\cO)$ a.s. for a subsequence $\{n_k\}$. It suffices to show that $I(X_{n})(T,\cO)=0$ a.s. on $A$ for all $n$.

{\em Step 4.} Assume that $X \in \cL_{\alpha}$ is bounded, continuous and satisfies (\ref{X-is-zero}). Let $(U_j^{(n)})_{j=1, \ldots,m_n}$ be a partition of $\cO$ in Borel sets with Lebesque measure smaller than $1/n$. Let $x_{j}^{n} \in U_j^{(n)}$ be arbitrary. Define
$$X_n(t,x)=\sum_{k=0}^{n-1}\sum_{j=1}^{m_n}X\left(\frac{kT}{n},x_{j}^{n}\right)1_{\left(\frac{kT}{n},
\frac{(k+1)T}{n}\right]}(t)
1_{U_j^{(n)}}(x).$$
Since $X$ is continuous in $(t,x)$, $X_n(t,x) \to X(t,x)$. By the bounded convergence theorem, $[X_n-X]_{\alpha} \to 0$, and hence $I(X_{n_k})(T,\cO) \to I(X)(T,\cO)$ a.s. for a subsequence $\{n_k\}$. Since on the event $A$,
$$I(X_n)(T,\cO)=\sum_{k=0}^{n-1} \sum_{j=1}^{m_n} X\left(\frac{kT}{n},x_{j}^{n}\right)
Z\left(\left(\frac{kT}{n},\frac{(k+1)T}{n}\right] \times U_j^{(n)}\right)=0,$$
it follows that $I(X)(T,\cO)=0$ a.s. on $A$.
$\Box$

\vspace{3mm}

\noindent \footnotesize{{\em Acknowledgement.} The author is grateful to Robert Dalang for suggesting this problem.

\normalsize{

\end{document}